\documentclass{jcmlatex}
\def\JCMvol{xx}
\def\JCMno{x}
\def\JCMyear{200x}
\def\JCMreceived{xxx}
\def\JCMrevised{xxx}
\def\JCMaccepted{xxx}

\usepackage{graphicx}
\usepackage{amsmath}
\usepackage{amsfonts}
\usepackage{graphics}
\usepackage{caption3}
\usepackage{amssymb}
\usepackage{psfrag}
\usepackage{amsmath}
\usepackage{graphicx}
\usepackage{multirow}
\usepackage{bm}
\usepackage{cite}
\usepackage{bm}
\usepackage{marginnote}
\usepackage{footnote}

\usepackage{graphicx}        % standard LaTeX graphics tool
\usepackage{multicol}        % used for the two-column index
\usepackage[bottom]{footmisc}% places footnotes at page bottom

\usepackage{amsmath,amssymb,bbm}
\usepackage{graphicx}
\usepackage[utf8]{inputenc}
\usepackage{url}
\usepackage{import}
\usepackage{color}
\usepackage{epic}
\usepackage{algpseudocode}
\usepackage{setspace}
\usepackage{tabularx}
\usepackage{array}

\newcommand{\ii}{\mathbf{i}}

\newtheorem{scheme}{Scheme}
\newtheorem{rem}{Remark}
\newtheorem{lem}{Lemma}

\begin{document}

\markboth{Weidong Zhao, Wei Zhang AND Guannan Zhang}{Template for Journal Of Computational
Mathematics}

\title{Second-order numerical schemes for decoupled forward-backward stochastic differential equations with jumps}

\author{Weidong Zhao
\thanks{Department of Mathematics, Shandong University, Jinan, Shandong 250100, China\\
Email: wdzhao@sdu.edu.cn}
\and Wei Zhang
\thanks{Department of Mathematics, Beijing University of Technology, Beijing 100022, China\\
Email:weizhang0313@bjut.edu.cn}
\and Guannan Zhang
\thanks{Department of Computational and Applied Mathematics, Oak Ridge National Laboratory, Oak Ridge, Tennessee 37831, United States\\ Email: zhangg@ornl.gov}
}

\maketitle

\begin{abstract}
We propose new numerical schemes for decoupled forward-backward stochastic differential equations (FBSDEs) with jumps, where the stochastic dynamics are driven by a $d$-dimensional Brownian motion and an independent compensated Poisson random measure.
A semi-discrete scheme is developed for discrete time approximation, which is constituted by a classic scheme for the forward SDE\cite{Platen:2010eo,Higham:2005gaa} and a novel scheme for the backward SDE.
Under some reasonable regularity conditions, we prove that the semi-discrete scheme can achieve second-order convergence in approximating the FBSDEs of interest; and such convergence rate does
not require jump-adapted temporal discretization.
Next, to add in spatial discretization, a fully discrete scheme is developed by designing accurate quadrature rules for estimating the involved conditional mathematical expectations.
Several numerical examples are given to illustrate the effectiveness and the high accuracy
 of the proposed schemes.
\end{abstract}



\section{Introduction}
In this work, we study numerical solution of decoupled forward-backward
stochastic differential equations (FBSDEs) with jumps,
where the underlying stochastic jump processes are characterized by Poisson random measures.
The term ``decoupled'' refers to the fact that the forward SDE
is independent of the solution of the backward SDE.
This work is motivated by a wide variety of applications offered by FBSDEs. In finance and insurance, FBSDEs-based
approaches \cite{Tang:1994gu,Peng:1993jt} have gained a great attention by both academics and practitioners, because FBSDEs provide us a unified framework to describe the mathematical problems which arise in option pricing \cite{ElKaroui:1997dn}, portfolio hedging \cite{EyraudLoisel:2005gk}, market utility maximization \cite{Becherer:2006jba} and risk measures \cite{Peng:2004ja,Quenez:2013il}, etc. Moreover, in the presence of jump behaviors in many financial problems \cite{Platen:2010eo}, L\`{e}vy jump processes have been incorporated into FBSDEs \cite{Delong:2013uk,EyraudLoisel:2005gk}, so as to accurately capture and properly interpret event-driven stochastic phenomena, such as corporate defaults, operational failures, insured events, etc.
In mathematics, one can relate FBSDEs with jumps to a class of nonlinear partial integro-differential equations (PIDEs), based on the extension of the nonlinear Feynman-Kac theory studied in \cite{Anonymous:fk}. As such, FBSDEs become a powerful probabilistic technique for studying analytical and numerical solutions and properties of the PIDEs, where the nonlocal integral operators of the PIDEs are characterized by Poisson random measures in the framework of FBSDEs. In engineering science, a particular application of the PIDEs is to model anomalous diffusion \cite{Metzler:2000du}, i.e.,~super-diffusion and sub-diffusion, that has been verified experimentally to be present in various applications, e.g.,~contaminant transport in groundwater and plasma physics. In this setting, FBSDEs-based probabilistic numerical schemes have been developed in \cite{BSDEnonlocal} to solve the governing PIDEs, which illustrated effectiveness of the FBSDEs models.

There are many theoretical results on FBSDEs with jumps over the past two decades.
The existence and uniqueness were proved by Tang and Li \cite{Tang:1994gu} for backward stochastic differential equations with Poisson jumps and Lipschitzian coefficients, which was then extended, by Rong in \cite{Rong:1997fe}, to the case of non-Lipschitzian coefficients.
In \cite{Anonymous:fk}, Barles, Buckdahn and Pardoux established a comparison theorem for decoupled FBSDEs with jumps as well as the link between such FBSDEs and PIDEs, which generalized the results in \cite{Pardoux:1990ju,Pardoux:1992jo} to the case of a natural filtration associated with a Brownian motion and a Poisson random measure. After that, in the context of FBSDEs with jumps, \O ksendal and Sulem  \cite{oksendal:2010kj} established maximum principles, and Royer \cite{Royer:2006ey} introduced nonlinear expectations. For a general overview of related topics, see \cite{Crepey:2013ki,Delong:2013uk} and the references therein.

The obstacle of applying FBSDEs with jumps to real-world engineering and finance problems results from the challenge of solving FBSDEs analytically or numerically. Since it is typically difficult to obtain analytical solutions, numerical solutions are highly desired in practical applications.
Numerical methods for FBSDEs without jumps have been well studied in the literature \cite{Douglas:1996em,Zhang:2004tt,ChenHong2013,Gobet:2005ft,Zhao:2014eu,
Zhao:2006jx,Zhao:2014hk,Bender:2007kz,Crisan:2014ed,Chassagneux:2014gn}, nevertheless, there are very few numerical schemes developed for FBSDEs with jumps, and most of those schemes only focused on temporal discretization. For instance, a Picard's iterative method was provided in \cite{Lejay:2007tu}, and numerical schemes of  backward SDE were studied in ~\cite{Bouchard:2008cp,Bouchard:2009il}.
%where the convergence rate was proved to be $(\Delta t)^{\frac{1}{2}}$. 
Due to the aforementioned applications of FBSDEs with jumps, it is of great significance to develop high-order temporal-spatial discretization schemes
for solving not only the FBSDEs but also the PIDEs and related  engineering problems.

In this paper, we propose novel numerical schemes for decoupled FBSDEs driven by a $d$-dimensional Brownian motion and an independent compensated Poisson random measure.
In general, the approximation of the FBSDEs under consideration includes two steps, i.e.,~constructing a semi-discrete scheme for temporal discretization, and extending it to a fully discrete scheme by incorporating effective spatial discretization.
%Specifically, in temporal discretization, the forward SDE can be discretized using classic numerical schemes (see \cite{Platen:2010eo,Higham:2005gaa} and the references therein) due to the decoupling property, and we develope a variant of the Crank-Nicolson scheme for the backward SDE.
By imposing appropriate regularity conditions on the coefficients, the generator  and the terminal condition, we rigorously prove the {\em second-order} convergence of the semi-discrete scheme with respect to $\Delta t$.
In spatial discretization, a carefully designed quadrature rule is critical to approximate all the involved conditional mathematical expectations which are, in this case, multiple integrals with respect to both the Brownian motion and the Poisson random measure. The integrals with respect to the Brownian motion is estimated by the Gauss-Hermite rule. For the integrals with respect to the Poisson random measure, we propose a general quadrature rule for the case that the jump component has {\em finite} activities.
A specific form of the quadrature rule can be determined based on the type of the underlying L\`{e}vy measure. For the numerical experiments in \S\ref{sec:num}, the L\`{e}vy measure is defined as a uniform distribution on bounded domains, so that Gauss-Legendre rule is an appropriate choice. Moreover, to avoid the explosion of the total number of quadrature points with the increase of time steps, we construct a piecewise Lagrange interpolating polynomials on a pre-determined spatial mesh, which are used to evaluate the integrand at all quadrature points.

%It is worth noting that the proposed schemes for the FBSDEs can be directly used to solve the terminal value problems for the related nonlinear PIDEs. Despite the probabilistic nature, our schemes are nevertheless {\em deterministic} algorithms. Thus, with the use of our approach for approximating the PIDEs, we can obtain the same convergence rate as that for the FBSDEs.
%%
%In fact,
%our method reveals several advantages over classic numerical approaches for PIDEs, e.g., finite element and finite difference methods. First, the proposed method includes a stable implicit time-stepping scheme, but do {\em not} require solving linear systems at each time step. This is signifiant in the sense of solving PIDEs, because the non-locality of the integral operators dramatically deteriorate the sparsity of the linear systems, which makes it extremely challenging to develop efficient linear solvers for the aforementioned classic numerical methods.
%%
%In addition, by using our scheme, the PIDE can be solved independently at different spatial grid points on each time level, making it straightforward to incorporate parallel implementation and adaptive spatial approximation.
The main contributions of this paper are as follows:
%\vspace{-0.05cm}
\begin{itemize}%\setlength{\itemsep}{0.00cm}
\item propose a second-order discrete time approximation (semi-discrete) scheme for decoupled FBSDEs with jumps.
\item rigorously analyze the convergence rate of the proposed discrete time approximation scheme with respect to $\Delta t$.
\item propose a fully discrete scheme by developing new quadrature rules for estimating involved conditional mathematical expectations.
%
%\item demonstrate that the fully discrete scheme is an effective and efficient technique for numerical solution of  the related PIDEs.
\end{itemize}

The outline of the paper is organized as follows. In \S\ref{sec:pre},
we introduce the mathematical description of the FBSDEs under consideration.
 In \S\ref{sec:ref}, we propose the semi-discrete scheme, i.e.,~time discretization, for the FBSDEs of interest.
Rigorous error analysis for the proposed semi-discrete scheme is conducted in \S\ref{sec:sta}.
The fully discrete scheme for the case of Poisson random measures with finite activities
is proposed in \S \ref{sec:fully}. Numerical examples are given in \S \ref{sec:num},
to show the effectiveness and the high accuracy of our approach.
Finally, several concluding remarks and discussions
about our future work are given in \S \ref{sec:con}.

\section{Preliminaries}\label{sec:pre}
Let $(\Omega, \mathcal{F}, \{\mathcal{F}_t\}_{0\le t \le T}, \mathbb{P})$
be a stochastic basis satisfying the usual hypotheses of completeness, i.e.,~$\mathcal{F}_0$
contains all the sets of $\mathbb{P}$-measure zero and possesses right continuity, i.e.,~$\mathcal{F}_t = \mathcal{F}_{t+}$. The filtration
$\{\mathcal{F}_t\}_{0\le t \le T}$ is assumed to be generated by two mutually independent processes,
i.e.,~one $d$-dimensional Brownian motion ${W}_t = (W_t^1, \ldots, W_t^d)^{\top}$
and one Poisson random measure $\mu(A,t)$ on $E \times [0,T]$ where
$E = \mathbb{R}^q \backslash\{0\}$ is equipped with its Borel field $\mathcal{E}$.
The compensator of $\mu$ and the resulting compensated Poisson random measure are
denoted by $\nu(de,dt) = \lambda(de)dt$ and $\tilde{\mu} (de, dt)= \mu(de, dt) - \lambda(de)dt$, respectively,
such that $\{\tilde{\mu}(A \times [0,t]) = (\mu - \nu)(A \times [0,t])\}_{0\le t \le T}$ is a martingale for all $A \in \mathcal{E}$.
$\lambda(de)$ is assumed to be a $\sigma$-finite measure on $(E, \mathcal{E})$ satisfying
\[
\int_{E} (1 \land |e |^2) \lambda (d e) < + \infty,
\]
where $|\cdot|$ denotes the standard Euclidean norm in Euclidean spaces.

In the probability space $(\Omega, \mathcal{F}, \{\mathcal{F}_t\}_{0\le t \le T}, \mathbb{P})$,
we introduce the following forward-backward stochastic differential equation with jumps
\begin{equation}\label{s1:e1}
\left\{\begin{aligned}
X_t & =   X_0+\int_0^t b(s,X_s)ds+\int_0^t \sigma(s,X_s)d W_s
        +\int_0^t \int_{E}c(s,X_{s-},e)\tilde{\mu}(de,ds),\\
Y_t & =   \xi+\int_t^T f(s,X_s,Y_s,Z_s,\Gamma_s) ds-\int_t^T Z_s dW_s
       -\int_t^T \int_{E}U_s(e)\tilde{\mu}(de,ds),
\end{aligned}\right.
\end{equation}
where the quadruplet $(X_t, Y_t, Z_t, U_t)$ is the unknown,
${b}: [0,T] \times \mathbb{R}^{q} \rightarrow \mathbb{R}^q$ is referred to as the drift coefficient,
${\sigma}: [0,T] \times \mathbb{R}^q \rightarrow \mathbb{R}^{q \times d}$ is referred to as the local diffusion coefficient,
${c}: [0,T] \times \mathbb{R}^q \times E \rightarrow \mathbb{R}^q$ is referred to as the jump coefficient,
${f}: [0,T] \times \mathbb{R}^q \times \mathbb{R}^{p} \times \mathbb{R}^{p\times d} \times \mathbb{R}^p \rightarrow \mathbb{R}^p$
is referred to as the generator of the FBSDE, and the process ${\Gamma}_s$
is defined by ${\Gamma}_s = \int_{E} {U}_s(e) \eta(e) \lambda(de)$
for a given bounded function $\eta: E \rightarrow \mathbb{R}$, i.e.,~$\sup_{e\in E}|\eta(e)| < +\infty$.
The terminal condition ${\xi}$ is an $\mathcal{F}_T$-measurable random vector in $\mathbb{R}^q$.
A quadruplet $(X_t, Y_t, Z_t, U_t)$ is called an $L^2$-adapted solution
if it is an $\{\mathcal{F}_t\}$-adapted, square integrable processes
satisfying the FBSDEs in \eqref{s1:e1}.

Under standard assumptions on the given data ${b}$, ${\sigma}$, ${f}$, ${\varphi}$ and $c$
(see \cite{Anonymous:fk} for details), there exists a unique solution
$(Y_t, Z_t, U_t) \in S^2 \times L^2(W) \times L^2(\tilde{\mu})$ for the backward SDE in \eqref{s1:e1},
where $S^2$ is the set of $\{\mathcal{F}_t\}$-adapted c\`{a}dl\`{a}g processes $\{Y_t, 0 \le t \le T\}$
such that
$
\|Y\|_{S^2}^2 := \mathbb{E}\left[\left( \sup_{0\le t \le T} |Y_t| \right)^2 \right]< \infty,
$
$L^2(W)$ the set of $\mathcal{F}_t$-progressively measurable $q\times d$ dimensional processes
$\{Z_t, 0 \le t \le T\}$ such that
$
\|Z\|_{L^2(W)}^2  := \mathbb{E} \left[\int_{0}^T |Z_t|^2 dt \right] < \infty,
$
and $L^2(\tilde{\mu})$ the set of mappings $U: \Omega \times [0,T] \times E \rightarrow \mathbb{R}$
such that
$
\| U\|_{L^2(\tilde{\mu})}^2 :=  \mathbb{E} \left[\int_0^T \int_E U_t(e)^2 \lambda(de) dt\right]  < \infty.
$

Now we introduce a class of {\em nonlinear} partial integro-differential equations (PIDEs)
that will be related to the FBSDEs in \eqref{s1:e1} later.
We consider the unique viscosity solution $u(t,x) \in \mathcal{C}([0,T]\times \mathbb{R}^q)$
of the following nonlinear PIDE, i.e.,
\begin{equation}\label{eq:Kol}
\left\{
\begin{aligned}
&\frac{\partial u}{\partial t}(t,x) + \widetilde{\mathcal{L}}[u](t,x) +{f}(t,x,u,{\sigma}\nabla u, \mathcal{B}[u])=0,
\text{ for } (t,x) \in [0,T) \times \mathbb{R}^q,\\
&u(T,x) = {\varphi}(x), \text{ for } x \in \mathbb{R}^q,
\end{aligned}
\right.
\end{equation}
where ${\varphi}(x)$ is the terminal condition at the time $t = T$, $\widetilde{\mathcal{L}}$ is the second-order integral-differential operator of the form
\begin{equation}\label{eq:L1}
\begin{aligned}
& \widetilde{\mathcal{L}}[u](t,x)  = \sum_{i=1}^q b_i(t,x)\frac{\partial u}{\partial x^i}(t,x)
+ \frac{1}{2}\sum_{i,j = 1}^q ({\sigma}{\sigma}^{\top})_{i,j}(t,x)
\frac{\partial^2 u}{\partial x^i \partial x^j}(t,x)\\
& \hspace{1.2cm} + \int_{E} \left(u(t,x + c(t,x,e)) -u(t,x) -
\sum_{i=1}^ q\frac{\partial u}{\partial x^i}(t,x) c(t,x,e) \right) \lambda (de),
\end{aligned}
\end{equation}
and $\mathcal{B}$ is an integral operator defined as
\[
\mathcal{B}[u](t,x) = \int_{E} \big[u(t,x + {c}(t,x,e)) -u(t,x)\big] \eta(e) \lambda(de).
\]
For $(t,x)\in [0,T]\times \mathbb{R}^q$, let $\mathbb{E}_{t}^x[\cdot]$ denote the mathematical expectation under the condition that
$X_t= x$, i.e.,~$\mathbb{E}_{t}^x[\cdot] := \mathbb{E}[\cdot | X_t = x]$.
To relate the FBSDEs in \eqref{s1:e1} with the PIDE in \eqref{eq:Kol},
we consider the FBSDEs of the following form
\begin{equation}\label{s3:e1}
\left\{\begin{aligned}
X_s^{t,x} &= x + \int_t^s b(r,X_r^{t,x})dr+\int_t^s \sigma(r,X_r^{t,x})d W_r
              + \int_t^s \int_E c(r,X_{r-}^{t,x},e) \tilde{\mu}(de,dr),\\
Y_s^{t,x} & =  \xi +\int_s^T f(r,X_r^{t,x},Y_r^{t,x},Z_r^{t,x},\Gamma_r^{t,x}) dr
               - \int_s^T Z_r^{t,x} dW_r
           -\int_s^T\int_E U_r^{t,x}(e)\tilde{\mu}(de,dr),
\end{aligned}\right.
\end{equation}
where the solution is $(X_s^{t,x},Y_s^{t,x},Z_s^{t,x},U_s^{t,x})$
and ${\Gamma}_s^{t,x} = \int_{E} {U}_s^{t,x}(e) \eta(e) \lambda(de)$ for $t\le s\le T$. Note that the superscripts
in \eqref{s3:e1} indicate the fact that the forward SDE in \eqref{s3:e1} starts from the
 time-space point $(t,x) \in [0,T] \times \mathbb{R}^q$.

According to Theorem 3.4 in \cite{Anonymous:fk}, if the terminal condition $\xi$
of the FBSDEs is a function of $X_T^{t,x}$, defined by $\xi = \varphi(X_T^{t,x})$ ($\varphi(\cdot)$
is the terminal condition of the PIDE), then the triple $(Y_s^{t,x}, Z_s^{t,x}, U_s^{t,x})$ for $t \le s \le T$
can be represented by the unique viscosity solution $u(t,x)$ of the PIDE \eqref{eq:Kol} as follows:
\begin{equation}\label{FK}
\left\{
\begin{aligned}
&Y_s^{t,x} = u(s,X_s^{t,x}),\\
&Z_s^{t,x} = {\sigma}(s, X_s^{t,x}) \nabla u(s, X_s^{t,x}), \\
& {U}_s^{t,x} = u(s,X_{s-}^{t,x}+{c}(s-,X_{s-}^{t,x},e)) - u(s, X_{s-}^{t,x}),
\end{aligned}\right.
\end{equation}
where $\nabla u$ denotes the gradient of $u$ with respect to $x$ and the function
$\Gamma_s^{t,x}$ is defined by $ {\Gamma}_s^{t,x} = \mathcal{B}[u](s, X_{s}^{t,x})$.
Particularly, when $s = t$, we have $u(t,x) = Y_t^{t,x} = \mathbb{E}[Y_t | X_t = x]$.

\section{The semi-discrete scheme for FBSDEs with jumps} \label{sec:ref}
In this section, we propose a numerical scheme for discrete-time
approximation of the FBSDEs under consideration.
Instead of the FBSDEs  \eqref{s1:e1}, we will use the conditional representation of the
FBSDEs given in \eqref{s3:e1} throughout this section. Specifically, discretizations of
the forward SDE and backward SDE are discussed in \S \ref{sec:DBSDE} and  \S \ref{sec:DFSDE}, respectively,
and the main numerical scheme is proposed in \S \ref{sec:sch}.
To proceed, we introduce the following time partition for the interval $[0,T]$:
\begin{equation}\label{tp}
\mathcal{T} := \{0=t_0<t_1<\cdots<t_{N}=T\}
\end{equation}
with $\Delta t_n := t_{n+1}-t_n$ and
$\Delta t := \max\limits_{0\leq n\leq N-1} \Delta t_n$.
We assume that the time partition $\mathcal{T}$ has the following regularity:
\begin{equation}\label{s3:e2}
\frac{\max\limits_{0\leq n\leq N-1} \Delta t_n}{\min\limits_{0\leq n\leq N-1} \Delta t_n}\leq c_0,
\end{equation}
where $c_0\geq 1$ is a real positive constant.
We remark  that $\mathcal{T}$ is {\em not} a jump-adapted partition.

\subsection{Discretization of the forward SDE}\label{sec:DBSDE}
Due to the decoupling of the FBSDEs in \eqref{s3:e1}, the forward SDE can be discretized separately.
Here we briefly recall some classic numerical schemes and their properties discussed in \cite{Platen:2010eo}.
Any of these schemes can serve as the approximation of the forward SDE in our schemes for the FBSDE.
By setting $t = t_n$, $s = t_{n+1}$ and $x = X^n$ in \eqref{s3:e1},
the forward SDE can be written as
\begin{equation}\label{XXX}
\begin{aligned}
X_{t_{n+1}}^{t_n,X^n}
= & X^n+\int_{t_n}^{t_{n+1}} b(s,X_s^{t_n,X^n}) ds
    + \int_{t_n}^{t_{n+1}} \sigma(s,X_s^{t_n,X^n}) d W_s\\
  & + \int_{t_n}^{t_{n+1}} \int_E c(s,X_{s-}^{t_n,X^n},e) \tilde{\mu}(de,ds),
\end{aligned}\end{equation}
where we assume that the solution $X_s^{t_n,X^n}$ starts at the time instant $t = t_n$
and spatial location $X_{t_n}^{t_n,X^n} = X^n$.
By using the It\^o-Taylor expansion, numerical schemes of strong-order $\beta$
(or the weak-order $\beta$) \cite{Platen:2010eo} can be represented in a general form, i.e.,
\begin{align}\label{s4:X111}
X^{n+1}
= & X^n + \Phi(t_{n},t_{n+1}, X^{n}, I_{\mathcal{J}\in \mathcal A_{\beta}}),
\end{align}
where $\Phi$ is the incremental, $\mathcal A_{\beta}$ is a hierarchical set such that the convergence
rate of the scheme is $\beta$ in a strong or weak sense.
Details of the index set $I_{\mathcal{J}\in \mathcal A_{\beta}}$ and
the definition of $\mathcal A_{\beta}$ can be found in \cite{Platen:2010eo} (pp.~196 and pp.~290).
The scheme \eqref{s4:X111} has the following properties:
\begin{itemize}
\item {\em Stability}: for an integer $r>0$,
there exists a constant $C\in (0,\infty)$ such that
\begin{equation}\label{s5:st}
\max\limits_{0\le n\le N}\mathbb{E}[|X^n|^r]\leq
C(1+\mathbb{E}[|X_0|^r]).
\end{equation}
\item {\em Approximation error:} there exist positive real numbers
$r_1, r_2, r_3, \alpha, \beta,\gamma$ such that for any function $g\in
\mathcal{C}_P^{2\beta+2}$, we have
\begin{equation}\label{s5:EX}\begin{aligned}
\left|\mathbb{E}^{X^n}_{t_n}[g(X_{t_{n+1}}^{t_n, X^n})-g(X^{n+1})]\right| & \leq
C(1+|X^n|^{2r_1})(\Delta t)^{\beta+1},\\
\left|\mathbb{E}^{X^n}_{t_n}[(g(X_{t_{n+1}}^{t_n, X^n})-g(X^{n+1}))\Delta
\tilde{W}^{\top}_{t_{n+1}}]\right| & \leq  C(1+|X^n|^{2r_2})(\Delta t)^{\gamma
+ 1},\\
\left|\mathbb{E}^{X^n}_{t_n}[(g(X_{t_{n+1}}^{t_n, X^n})-g(X^{n+1}))
\Delta\tilde{\mu}_{t_{n+1}}^{*}]\right|
&\leq  C(1+|X^n|^{2r_3})(\Delta t)^{\alpha+ 1},
\end{aligned}\end{equation}
where $C_P^{2\beta+2}$ is the set of $2\beta+2$ times continuously differentiable functions which,
together with their derivatives of order up to
$2\beta+2$, have at most polynomial growth. According to Theorem 6.4.1
%on page 201
and Theorem 12.3.4
%on page 364
in \cite{Platen:2010eo}, it is easy to derive that $\alpha = \beta = \gamma$ for strong and weak Taylor schemes.
In this paper, we prove in Theorem \ref{s5:TL2} that the second-order convergence of the proposed semi-discrete scheme for the FBSDE
in \eqref{s1:e1} requires $\alpha=\beta=\gamma=2$.
\end{itemize}

\subsection{Discretization of the backward SDE}\label{sec:DFSDE}

Now we study the discretization of the backward SDE in \eqref{s3:e1}
driven by the process $X_{s}^{t_n, X^n}$ in \eqref{XXX} for
$s \in [t_n, t_{n+1}]$. Within the interval $[t_n, t_{n+1}]$, the backward SDE can be rewritten as
\begin{equation}\label{BSDEtn}
\begin{aligned}
Y_{t_n}^{t_n, X^n}
=  Y_{t_{n+1}}^{t_n, X^n}+\int_{t_n}^{t_{n+1}} f_s^{t_n, X^n}ds
-\int_{t_n}^{t_{n+1}} Z_s^{t_n, X^n} dW_s
- \int_{t_n}^{t_{n+1}} \int_E U_s^{t_n, X^n}(e)\tilde{\mu}(de,ds),
\end{aligned}
\end{equation}
where $f_s^{t_n, X^n}$ denotes $f(s,X_s^{t_n, X^n},Y_s^{t_n, X^n},Z_s^{t_n, X^n},\Gamma_s^{t_n, X^n})$
for notational simplicity.
Due to the relation between $\Gamma_s^{t,x}$ and $U_{s}^{t,x}$, in what follows, all the numerical schemes for the backward SDE will be proposed to
approximate $(Y_s^{t_n,X^n}, Z_s^{t_n,X^n}, \Gamma_s^{t_n,X^n})$.
Since there are three unknown stochastic processes involved in \eqref{BSDEtn},
we now construct three discretized reference equations for $Y_{t_n}^{t_n, X^n}$,
$Z_{t_n}^{t_n, X^n}$ and $\Gamma_{t_n}^{t_n, X^n}$ in \S \ref{refY}, \S \ref{refZ}
and \S \ref{refG}, respectively, which are the foundation of the formal
semi-discrete scheme discussed in \S \ref{sec:sch} for the FBSDEs.

\subsubsection{The reference equation for $Y_{t_n}^{t_n,X^n}$}\label{refY}
Taking the conditional mathematical expectation
$\mathbb{E}_{t_n}^{X^n} [\cdot]$ on both sides  of \eqref{BSDEtn}, we obtain
\begin{equation}\label{s3:Y1:1}
Y_{t_n}^{t_n, X^n}=\mathbb{E}_{t_n}^{X^n}\left[Y_{t_{n+1}}^{t_n, X^n}\right]
+\int_{t_n}^{t_{n+1}} \mathbb{E}_{t_n}^{X^n}\left[f_s^{t_n, X^n}\right]ds,
\end{equation}
due to the fact that $\int_{t_n}^{t_{n+1}} Z_s^{t_n, X^n} dW_s$
and $\int_{t_n}^{t_{n+1}} \int_E U_s^{t_n, X^n}(e) \tilde{\mu}(de,ds)$ for $t > t_n$ are martingales.
Note that the integrand
$\mathbb{E}_{t_n}^{X^n}[f_s^{t_n, X^n}]$
is a {\em deterministic} function of $s \in [t_n, t_{n+1}]$ under the $\sigma$-algebra $\mathcal{F}_{t_n}$. Thus,
numerical integration approaches can be used to approximate the temporal integral in \eqref{s3:Y1:1}.
In this effort, we use the Crank-Nicolson scheme, i.e.,~the trapezoidal rule, such that
\begin{equation}\label{s3:e5}
\int_{t_n}^{t_{n+1}}\mathbb{E}_{t_n}^{X^n}\left[f_s^{t_n, X^n}\right] ds
=\frac{1}{2}\Delta t_n f_{t_n}^{t_n, X^n}+\frac{1}{2}\Delta t_n \mathbb{E}_{t_n}^{X^n}\left[f_{t_{n+1}}^{t_n, X^n}\right]
+R_y^n,
\end{equation}
where the residual $R_y^n$ is
\begin{equation}\label{s3:Ry}
R_y^n :=\int_{t_n}^{t_{n+1}}\bigg\{\mathbb{E}_{t_n}^{X^n}\left[f_s^{t_n,X^n}\right]-\frac{1}{2}
f_{t_n}^{t_n,X^n}-\frac{1}{2} \mathbb{E}_{t_n}^{X^n}\left[f_{t_{n+1}}^{t_n,X^n}\right]\bigg\}ds.
\end{equation}
Substituting (\ref{s3:e5}) into (\ref{s3:Y1:1}), we obtain the reference equation for solving $Y_{t_n}^{t_n,X^n}$:
\begin{equation}\label{s3:e6}
Y_{t_n}^{t_n,X^n}=\mathbb{E}_{t_n}^{X^n}\left[Y_{t_{n+1}}^{t_n,X^n}\right]
+\frac{1}{2}\Delta t_n f_{t_n}^{t_n,X^n}+\frac{1}{2}\Delta t_n \mathbb{E}_{t_n}^{X^n}\left[f_{t_{n+1}}^{t_n,X^n}\right]
+R_y^n.
\end{equation}

\subsubsection{The reference equation for $Z_{t_n}^{t_n, X^n}$}\label{refZ}
To proceed, we introduce a new Gaussian process $\Delta \tilde {W}_s$ defined by
\begin{equation}\label{s3:NBM}
\Delta \tilde{W}_s=2\Delta W_s-\frac{3}{\Delta t_n}\int_{t_n}^{s}(r-t_n)dW_r, \quad \forall s \in [t_n, t_{n+1}],
\end{equation}
where $\Delta W_s=W_s-W_{t_n}$ is the $d$-dimensional standard Brownian motion in the FBSDEs in \eqref{s3:e1}.
It is easy to see that
$\Delta \tilde{W}_s=(\Delta \tilde{W}_s^1,\Delta \tilde{W}_s^2,\cdots,\Delta \tilde{W}_s^d)^{\top}$
is also a $d$-dimensional Gaussian process with the properties
$\mathbb{E}_{t_n}^{X^n}[\Delta \tilde{W}_s]=0$, $
\mathbb{E}_{t_n}^{X^n}[\Delta \tilde{W}_s^i\Delta \tilde{W}_s^j]=0$ for $i\ne j$, and
\begin{equation*}
\begin{aligned}
\mathbb{E}_{t_n}^{X^n}[(\Delta \tilde{W}_s^i)^2]
&=\mathbb{E}_{t_n}^{X^n}\Big[(2\Delta W_s^i-\frac{3}{\Delta t_n}\int_{t_n}^s(r-t_n)dW_r^i)^2\Big]\\
&=4(s-t_n)-\frac{6(s- t_n)^2}{\Delta t_n}+\frac{3(s- t_n)^3}{
\Delta t_n^2}, \;\; \text{ for } i = 1, \ldots, d.
\end{aligned}
\end{equation*}
In the case of $s = t_{n+1}$, we have $\mathbb{E}_{t_n}^{X^n}[\Delta \tilde {W}^i_{t_{n+1}}]=0$ and
$\mathbb{E}_{t_n}^{X^n}[(\Delta \tilde{W}^i_{t_{n+1}})^2]=\Delta t_n$ for $i = 1, \ldots, d$.

Multiplying $\eqref{BSDEtn}$ by the transpose of $\Delta \tilde{W}_{t_{n+1}}$ in \eqref{s3:NBM},
and taking the conditional mathematical expectation
$\mathbb{E}_{t_n}^{X^n}[\cdot]$ on both sides, we obtain
\begin{equation}\label{s3:e7}
\begin{aligned}
0 & = \mathbb{E}_{t_n}^{X^n}\left[Y_{t_{n+1}}^{t_n,X^n}\Delta \tilde{W}_{t_{n+1}}^{\top}\right]
+\int_{t_n}^{t_{n+1}} \mathbb{E}_{t_n}^{X^n}\left[f_s^{t_n,X^n} \Delta \tilde{W}_{t_{n+1}}^{\top}\right]ds\\
&\quad -\mathbb{E}_{t_n}^{X^n}\left[\int_{t_n}^{t_{n+1}}Z_s^{t_n,X^n} dW_s\cdot\Delta \tilde{W}_{t_{n+1}}^{\top}\right].
\end{aligned}
\end{equation}
Then, the right endpoint rule is used to discretize the first temporal integral in \eqref{s3:e7}, such that
\begin{equation}\label{s3:e8}
\begin{aligned}
\int_{t_n}^{t_{n+1}} \mathbb{E}_{t_n}^{X^n}\left[f_s^{t_n,X^n} \Delta \tilde{W}_{t_{n+1}}^{\top}\right]ds
&= \Delta t_n \mathbb{E}_{t_n}^{X^n}\left[f_{t_{n+1}}^{t_n,X^n}\Delta \tilde{W}_{t_{n+1}}^{\top}\right]
   + R_{z,1}^n,\\
\end{aligned}
\end{equation}
where $R_{z,1}^n:=\int_{t_n}^{t_{n+1}}
\mathbb{E}_{t_n}^{X^n}[f_s^{t_n,X^n} \Delta \tilde{W}_{t_{n+1}}^{\top}]ds
-\Delta t_n \mathbb{E}_{t_n}^{X^n}[f_{t_{n+1}}^{t_n,X^n}\Delta \tilde{W}_{t_{n+1}}^{\top}]$ is the residual.
For the second temporal integral in \eqref{s3:e7},
based on the properties of $\Delta \tilde{W}^{\top}_{t_{n+1}}$,
we discretize it by
\begin{equation}\label{s3:e9}
\begin{aligned}
-\mathbb{E}_{t_n}^{X^n}\left[\int_{t_n}^{t_{n+1}}Z_s^{t_n,X^n} dW_s\cdot\Delta \tilde{W}_{t_{n+1}}^{\top}\right]
%= & -Z_{t_n}^{t_n,X^n}\mathbb{E}^{X^n}_{t_n}\left[\Delta W_{t_{n+1}}\Delta \tilde W_{t_{n+1}}^{\top}\right]+R_{z,2}^n
= -\frac{1}{2}\Delta t_n Z_{t_n}^{t_n,X^n}+R_{z,2}^n,
\end{aligned}
\end{equation}
where the residual is $ R_{z,2}^n:=\frac{1}{2}\Delta t_n Z_{t_n}^{t_n,X^n}
-\mathbb{E}_{t_n}^{X^n}[\int_{t_n}^{t_{n+1}}Z_s^{t_n,X^n} dW_s\cdot\Delta \tilde{W}_{t_{n+1}}^{\top}].$
Substituting (\ref{s3:e8}) and (\ref{s3:e9}) into (\ref{s3:e7}),
we obtain the reference equation for $Z_{t_n}^{t_n,X^n}$, i.e.,
\begin{equation}\label{s3:e10}
\frac{1}{2}\Delta t_n Z_{t_n}^{t_n,X^n}=\mathbb{E}_{t_n}^{X^n}
\left[Y_{t_{n+1}}^{t_n,X^n}\Delta \tilde{W}_{t_{n+1}}^{\top}\right]
+\Delta t_n  \mathbb{E}_{t_n}^{X^n}\left[f_{t_{n+1}}^{t_n,X^n}\Delta \tilde{W}_{t_{n+1}}^{\top}\right]+R_{z}^n,
\end{equation}
where $R_{z}^n:=R_{z,1}^n+R_{z,2}^n$.

\subsubsection{The reference equation for $\Gamma_{t_n}^{t_n, X^n}$} \label{refG}
Similar to the definition of $\Delta \tilde W_s$,
by using the compensated Poisson random measure $\tilde{\mu}(de,ds)$ in \eqref{s3:e1},
we define a new stochastic process $\Delta\tilde{\mu}_{s}^*$ as
\begin{equation}\label{s3:Np}
\Delta\tilde{\mu}_{s}^*
=\int_{t_n}^s\int_E \left(2-\frac{3(t-t_n)}{\Delta t_n}\right)\eta(e)\tilde{\mu}(de,dt),
\quad \forall s \in [t_n, t_{n+1}].
\end{equation}
Then, multiplying \eqref{BSDEtn}
by $\Delta\tilde{\mu}_{t_{n+1}}^{*}$ and
taking the conditional mathematical expectation
$\mathbb{E}_{t_n}^{X^n}[\cdot]$ on both sides, we obtain
\begin{equation}\label{s3:e11}
\begin{aligned}
0 & =\mathbb{E}_{t_n}^{X^n}\left[Y_{t_{n+1}}^{t_n,X^n}\Delta\tilde{\mu}_{t_{n+1}}^{*}\right]
    +\int_{t_n}^{t_{n+1}}\mathbb{E}_{t_n}^{X^n}\left[f_s^{t_n,X^n}\Delta\tilde{\mu}_{t_{n+1}}^{*}\right]ds\\
  & \hspace{0.4cm} -\mathbb{E}_{t_n}^{X^n}\left[\int_{t_n}^{t_{n+1}}
       \int_E U_s^{t_n,X^n}(e)\tilde{\mu}(de,ds)\Delta\tilde{\mu}_{t_{n+1}}^{*}\right].
\end{aligned}
\end{equation}
Analogous to the reference equation \eqref{s3:e10},
we also discretize the first temporal integral in \eqref{s3:e11} using the right endpoint rule, such that
\begin{equation}\label{s3:e12}
\int_{t_n}^{t_{n+1}}\mathbb{E}_{t_n}^{X^n}
\left[f_s^{t_n,X^n}\Delta\tilde{\mu}_{t_{n+1}}^{*}\right]ds
=\Delta t_n \mathbb{E}_{t_n}^{X^n}
\left[f_{t_{n+1}}^{t_n,X^n}\Delta\tilde{\mu}_{t_{n+1}}^{*}\right]+R_{\Gamma, 1}^n,
\end{equation}
where
$R_{\Gamma,1}^n:= \int_{t_n}^{t_{n+1}}\mathbb{E}_{t_n}^{X^n}[f_s^{t_n,X^n}\Delta\tilde{\mu}_{t_{n+1}}^{*}]ds
 -\Delta t_n \mathbb{E}_{t_n}^{X^n}[f_{t_{n+1}}^{t_n,X^n}\Delta\tilde{\mu}_{t_{n+1}}^{*}]
$ is the residual. For the second temporal integral in \eqref{s3:e11}, we have
\begin{equation}\label{s3:e13}
\begin{aligned}
&-\mathbb{E}_{t_n}^{X^n}\left[\int_{t_n}^{t_{n+1}}\int_E U_s^{t_n,X^n}(e)\tilde{\mu}(de,ds)
\Delta\tilde{\mu}_{t_{n+1}}^{*}\right]\\
= & -\mathbb{E}_{t_n}^{X^n}\left[\int_{t_n}^{t_{n+1}}\int_E U_{t_n}^{t_n,X^n}(e)\tilde{\mu}(de,ds)
\Delta\tilde{\mu}_{t_{n+1}}^{*}\right]+R_{\Gamma,2}^n\\
= & -\mathbb{E}_{t_n}^{X^n}\left[\Gamma_{t_n}^{t_n,X^n}
\int_{t_n}^{t_{n+1}}\left(2-\frac{3(s-t_n)}{\Delta t_n}\right)ds\right]+R_{\Gamma,2}^n\\
= & -\frac{1}{2}\Delta t_n \Gamma_{t_n}^{t_n,X^n}+R_{\Gamma, 2}^n,
\end{aligned}
\end{equation}
where
$R_{\Gamma, 2}^n
:= \frac{1}{2}\Delta t_n \Gamma_{t_n}^{t_n,X^n}
  - \mathbb{E}_{t_n}^{X^n}[\int_{t_n}^{t_{n+1}}\int_E U_s^{t_n,X^n}(e)\tilde{\mu}(de,ds)
\Delta\tilde{\mu}_{t_{n+1}}^{*}].$
By (\ref{s3:e11}), (\ref{s3:e12}) and (\ref{s3:e13}),
we obtain the reference equation for $\Gamma_{t_n}^{t_n,X^n}$, i.e.,
\begin{equation}\label{s3:e14}
\begin{aligned}
\frac{1}{2}\Delta t_n \Gamma_{t_n}^{t_n,X^n}
=\mathbb{E}_{t_n}^{X^n}\left[Y_{t_{n+1}}^{t_n,X^n}\Delta\tilde{\mu}_{t_{n+1}}^{*}\right]
 +\Delta t_n \mathbb{E}_{t_n}^{X^n}\left[f_{t_{n+1}}^{t_n,X^n}\Delta\tilde{\mu}_{t_{n+1}}^{*}\right]
 +R_{\Gamma}^n,
\end{aligned}
\end{equation}
where $R_{\Gamma}^n:=R_{\Gamma,1}^n+R_{\Gamma,2}^n$.

\subsection{The semi-discrete scheme}\label{sec:sch}
Now we combine the approximation in \eqref{s4:X111} of the forward SDE and the three reference equations in  \eqref{s3:e6}, \eqref{s3:e10} and \eqref{s3:e14} to propose our semi-discrete scheme (temporal discretization) for the FBSDEs in \eqref{s1:e1}.
Let $(X^{n+1}, Y^n, Z^n,\Gamma^n)$ denote the approximation to the exact solution
$(X_{t_{n+1}}^{t_n, X^n} Y_{t_n}^{t_n,X^n},Z_{t_n}^{t_n,X^n},\Gamma_{t_n}^{t_n,X^n})$ of the FBSDEs in \eqref{s1:e1} for $n=N-1,\ldots,0$. Based on the partition $\mathcal{T}$ of the time interval $[0,T]$, the approximate solution $(X^{n+1}, Y^n, Z^n,\Gamma^n)$ is constructed following the procedure in Scheme \ref{s4:sch:1}.

\begin{scheme}\label{s4:sch:1}
Given the initial condition $X_0$ for the forward SDE in \eqref{s1:e1}
and the terminal condition $(Y^N, Z^N,\Gamma^N)$ for the backward SDE in \eqref{s1:e1}, solve the approximate solution
$(X^{n+1}, Y^n, Z^{n},\Gamma^n)$, for $n = N-1, \ldots, 0$, by
\begin{align}
\label{s4:X}
X^{n+1} & =  X^n + \Phi(t_{n},t_{n+1}, X^{n}, I_{\mathcal{J}\in \mathcal A_{\beta}}),\\
\label{s4:Z}
\frac{1}{2}\Delta t_n Z^n
& =  \mathbb{E}_{t_n}^{X^n}\left[Y^{n+1}\Delta \tilde{W}_{t_{n+1}}^{\top}\right]
    +\Delta t_n \mathbb{E}_{t_n}^{X^n}\left[f^{n+1}\Delta \tilde{W}_{t_{n+1}}^{\top}\right],\\
\label{s4:G}
\frac{1}{2}\Delta t_n \Gamma^n
& = \mathbb{E}_{t_n}^{X^n}\left[Y^{n+1}\Delta\tilde{\mu}_{t_{n+1}}^{*}\right]
    +\Delta t_n \mathbb{E}_{t_n}^{X^n}\left[f^{n+1}\Delta\tilde{\mu}_{t_{n+1}}^{*}\right],\\
\label{s4:Y}
Y^n & = \mathbb{E}_{t_n}^{X^n}\left[Y^{n+1}\right]+\frac{1}{2}\Delta t_n f^n
    + \frac{1}{2}\Delta t_n \mathbb{E}_{t_n}^{X^n}\left[f^{n+1}\right],
\end{align}
where $f^{n+1} = f(t_{n+1},X^{n+1}, Y^{n+1}, Z^{n+1},
\Gamma^{n+1})$, $f^n = (t_n, X^n, Y^n, Z^{n}, \Gamma^{n})$,
$\Delta \tilde W_{t_{n+1}}$ and $\Delta\tilde{\mu}_{t_{n+1}}^{*}$ are defined according to
\eqref{s3:NBM} and \eqref{s3:Np} by setting $s = t_{n+1}$, respectively.
\end{scheme}

From the dependence of \eqref{s4:X}--\eqref{s4:Y}, we can see
that the scheme in \eqref{s4:X} for $X^{n+1}$ is independent of the other three schemes,
so that, at each time step, $X^{n+1}$ is always firstly determined.
Then, by observing that \eqref{s4:Z} and \eqref{s4:G} are {\em explicit} schemes,
we can solve $Z^n$ and $\Gamma^n$ by substituting $X^{n+1}$ into \eqref{s4:Z} and \eqref{s4:G}, respectively.
Next, since \eqref{s4:Y} includes $f^n$ that depends on $Y^n$, $Z^n$ and $\Gamma^n$,
it is an {\em implicit} scheme for $Y^n$. If the generator $f^n$ is {\em nonlinear}
and Lipschitz continuous with respect to $Y^n$, then $Y^n$ can be obtained
by substituting $X^{n+1}$, $Z^n$ and $\Gamma^n$ into \eqref{s4:Y} and solving a nonlinear equation.

In addition, we would like to discuss the application of Scheme \ref{s4:sch:1}
to discrete time approximations of the PIDE in \eqref{eq:Kol} when the terminal condition $\xi$
of the FBSDEs is a function of $X_T$, i.e.,~$\xi = \varphi(X_T)$.
The goal is to construct an approximate solution $u(t_n, x)$ for $n = N-1, \ldots, 1, 0$ and $x \in \mathbb{R}^q$.
Specifically, based on the relation between
$u(t,x)$ and $(X_t, Y_t, Z_t, \Gamma_t)$ in \eqref{FK}, the discrete time approximation, denoted by $u^n(x)$, is defined by
\begin{equation}\label{un}
u^n(x) := \mathbb{E}\left[Y^n | X^n = x\right] \approx u(t_n, x)
= \mathbb{E} \left[Y_{t_n} | X_{t_n} = x\right] \;\; \text{ for } \; x \in \mathbb{R}^q.
\end{equation}
It is easy to see that both $Y_{t_n}$ and $Y^n$ are {\em deterministic} values
under the conditions $X_{t_n} = x$ and $X^n = x$, respectively.
Moreover, the convergence of $(X^{n+1}, Y^{n}, Z^{n}, \Gamma^n)$ to $(X_{t_{n+1}}, Y_{t_n}, Z_{t_n}, \Gamma_{t_n})$
as $N \rightarrow \infty$ will ensure the convergence of $u^n(x)$ to $u(t_n, x)$.
Hence, Scheme \ref{s4:sch:1} can be viewed as an effective probabilistic scheme for the PIDE in \eqref{eq:Kol}.
Moreover, $Z^n$ and $\Gamma^n$ provide approximations of $\sigma \nabla u$ and  $\mathcal{B}[u]$ which enable accurate characterization of local and nonlocal diffusive fluxes in practical engineering problems.
%On the other hand, if only $u(t,x)$ is desired and $f$ is independent of $Z^n$ and $\Gamma^n$,
%then $u^n(x)$ can be obtained using only \eqref{s4:X} and \eqref{s4:Y} without solving $Z^n$ and $\Gamma^n$.

\section{Error estimates for the semi-discrete scheme}\label{sec:sta}
In this section, we estimate the truncation error of Scheme \ref{s4:sch:1}.
Since error estimates for the scheme for the forward SDEs have been well established in the literature
(see \cite{Platen:2010eo} and the references therein), we focus on analyzing the approximation error
of $(Y^n, Z^n, \Gamma^n)$ for $n = 0, \ldots, N-1$.
The general procedure of our analysis is similar to that for classic time-stepping schemes.
We first construct an upper bound of the global truncation error of $(Y^n, Z^n, \Gamma^n)$
by recursively accumulating local truncation errors.
Then, we estimate all the local truncation errors in the upper bound,
which relates the global truncation error to the maximum time step size
$\Delta t := \max\limits_{1\le n \le N} \Delta t_n$.

To proceed, we need to specify the definition of the approximation error of $(Y^i, Z^i, \Gamma^i)$.
For $i = 0, \ldots, N$, the errors of $Y^i$, $Z^i$, $\Gamma^i$
and $f^i := f(t_i, X^i, Y^i, Z^i, \Gamma^i)$ are respectively defined by
\begin{equation}\label{err}
\begin{aligned}
e_y^i := {Y}_{t_{i}}^{t_i,X^i}-Y^{i}, & \;\;\;e_z^i := {Z}_{t_{i}}^{t_i,X^i}-Z^{i}, & \\
e_\Gamma^i := {\Gamma}_{t_{i}}^{t_i,X^i}-\Gamma^{i}, & \;\;\;e_f^i := {f}_{t_{i}}^{t_i,X^i}-f^{i}, &\\
\end{aligned}
\end{equation}
where $Y_{t_i}^{t_i,X^i} := \mathbb{E}[ Y_{t_i} | X_{t_i} = X^i]$
and likewise for $Z_{t_i}^{t_i,X^i}$, $\Gamma_{t_i}^{t_i,X^i}$, $f_{t_i}^{t_i,X^i}$.
It should be noted that $Y_{t_{i+1}}^{t_i, X^i}$ and $Y_{t_{i+1}}^{t_{i+1}, X^{i+1}}$ for $0 \le i \le N$ are usually different stochastic processes because of the difference between $X_{t_{i+1}}^{t_i, X^i}$ and $X^{i+1}$. This fact can be easily shown with the use of the solution $u(t,x)$ of the PIDE in \eqref{eq:Kol}. According to the relationship in \eqref{FK}, it is easy to see that
$$Y_{t_{i+1}}^{t_i,X^i}=u(t_{i+1},X_{t_{i+1}}^{t_i,X^i}), \quad
Y_{t_{i+1}}^{t_{i+1},X^{i+1}}=u(t_{i+1},X^{i+1}),$$
where $X_{t_{i+1}}^{t_i,X^i}$ and $X^{i+1}$ are obtained by \eqref{XXX} and \eqref{s4:X}, respectively. As such, we introduce the following residual notations that will be used later:
\begin{equation}\label{s5:TL:e1aa}\begin{aligned}
\mathcal{R}_{y_1}^i & := \mathbb{E}_{t_i}^{X^i}\left[Y_{t_{i+1}}^{t_i,X^i}-{Y}_{t_{i+1}}^{t_{i+1},X^{i+1}}\right],
\\
\mathcal{R}_{y_2}^i  &:= \mathbb{E}_{t_i}^{X^i}\left[f_{t_{i+1}}^{t_i,X^i}- {f}_{t_{i+1}}^{t_{i+1},X^{i+1}}\right],\\
\mathcal{R}_{z_1}^i & := \mathbb{E}_{t_i}^{X^i}\left[\left(Y_{t_{i+1}}^{t_i,X^i}-{Y}_{t_{i+1}}^{t_{i+1},X^{i+1}}\right)
              \Delta \tilde{W}_{t_{i+1}}^{\top}\right],\\
\mathcal{R}_{z_2}^i & := \mathbb{E}_{t_i}^{X^i}\left[\left(f_{t_{i+1}}^{t_i,X^i}-{f}_{t_{i+1}}^{t_{i+1},X^{i+1}}\right)
              \Delta \tilde{W}_{t_{i+1}}^{\top}\right],\\
\mathcal{R}_{\Gamma_1}^i & := \mathbb{E}_{t_i}^{X^i}\left[\left(Y_{t_{i+1}}^{t_i,X^i}
                   - {Y}_{t_{i+1}}^{t_{i+1},X^{i+1}}\right)\Delta\tilde{\mu}_{t_{i+1}}^{*}\right],
\\
\mathcal{R}_{\Gamma_2}^i & := \mathbb{E}_{t_i}^{X^i}\left[\left(f_{t_{i+1}}^{t_i,X^i}
                   - {f}_{t_{i+1}}^{t_{i+1},X^{i+1}}\right)\Delta\tilde{\mu}_{t_{i+1}}^{*}\right]
\end{aligned}\end{equation}
for $i = 0, \ldots, N-1$. Note that the above residuals represent the local weak approximations
of the scheme \eqref{s4:X} for solving forward SDEs.

In the following theorem, we construct an upper bound of the errors
$e_y^n$, $e_z^n$ and $e_\Gamma^n$ $(0 \le n \le N-1)$
with the use of the residuals $R_y^i$, $R_z^i$ and $R_{\Gamma}^i$
in \eqref{s3:e6}, \eqref{s3:e10} and \eqref{s3:e14}, respectively, for $i = n, \ldots, N$,
as well as the residuals defined in \eqref{s5:TL:e1aa}.

\begin{theorem}\label{s5:TL}
Based on the partition $\mathcal{T}$ in \eqref{tp} of the time interval $[0,T]$,
if the generator $f(t,x,y,z,\gamma)$ is Lipschitz continuous with respect to
$x$, $y$, $z$ and $\gamma$ where Lipschitz constant is denoted by $L$,
then with sufficiently small time step $\Delta t : = \max\limits_{1\le i \le N} \Delta t_i$,
the errors $e_y^n$, $e_z^n$ and $e_\Gamma^n$ in \eqref{err} for $n = 0, \ldots, N-1$ can be bounded by
{\allowdisplaybreaks
\begin{equation}\label{s5:TL:e1}
\begin{aligned}
&\mathbb{E}[|e_y^n|^2]+ \Delta t\sum\limits_{i=n}^{N-1}
                        \bigg(\frac{1+C\Delta t}{1-C\Delta t}\bigg)^{i-n}\mathbb{E}[|e_z^i|^2 + |e_\Gamma^i|^2]\\
 \leq & \;C' \Big(\mathbb{E}[|e_y^N|^2] + \Delta t \mathbb{E}[|e_z^N|^2+|e_{\Gamma}^N|^2]\Big)\\
& + \sum\limits_{i=n}^{N-1}\bigg(\frac{1+C\Delta t}{1-C\Delta t}\bigg)^{i-n}\frac{C}{1-C\Delta t}
            \mathbb{E}\Bigg[
           \frac{1}{\Delta t} \left(|R_y^i|^2 + |R_{z}^i|^2 + |R_{\Gamma}^i|^2 \right)\\
& +\frac{1}{\Delta t}\left(|\mathcal{R}_{y_1}^i|^2+|\mathcal{R}_{z_1}^i|^2+ |\mathcal{R}_{\Gamma_1}^i|^2\right)
          + \Delta t\left( |\mathcal{R}_{y_2}^i|^2+ |\mathcal{R}_{z_2}^i|^2 + |\mathcal{R}_{\Gamma_2}^i|^2\right)\Bigg],
\end{aligned}
\end{equation}
}
where $C$ is a positive constant depending on $L$ and $c_0$ defined in \eqref{s3:e2}, $C'$
is a positive constant depending on $c_0$, $T$ and $L$, the residuals $R_y^i$, $R_{z}^i$
and $R_\Gamma^i$ for $i = n, \ldots, N$ are defined in \eqref{s3:e6}, \eqref{s3:e10} and \eqref{s3:e14}, respectively,
and $\mathcal{R}_{y_1}^i$, $\mathcal{R}_{y_2}^i$, $\mathcal{R}_{z_1}^i$, $\mathcal{R}_{z_2}^i$, $\mathcal{R}_{\Gamma_1}^i$, $\mathcal{R}_{\Gamma_2}^i$ are defined in \eqref{s5:TL:e1aa}.
\end{theorem}
\begin{proof}
This proof consists of four steps. Step 1, 2 and 3 are dedicated to the estimating
$e_y^n$, $e_z^n$ and $e_\Gamma^n$, respectively, and
those estimates are combined together at Step 4 that completes the proof.

$\bullet$ {\bf  Step 1}:~{\em Estimating the error $e_y^n = {Y}_{t_{n}}^{t_n,X^n}-Y^{n}$.}

Subtracting the scheme  \eqref{s4:Y} from the reference equation \eqref{s3:e6}, we have
\begin{equation*}\label{s5:TL:e2}
\begin{aligned}
e_y^n= \; & \mathbb{E}_{t_n}^{X^n}\Big[Y_{t_{n+1}}^{t_n,X^n}-Y^{n+1}\Big]
+ \frac{1}{2}\Delta t_n\left(f_{t_n}^{t_n,X^n}-f^n\right)\\
& + \frac{1}{2}\Delta t_n\mathbb{E}_{t_n}^{X^n}\Big[f_{t_{n+1}}^{t_n,X^n}-f^{n+1}\Big]+ R_y^n\\
=\;&\mathbb{E}_{t_n}^{X^n}\Big[Y_{t_{n+1}}^{t_n,X^n}-Y_{t_{n+1}}^{t_{n+1},X^{n+1}}
+Y_{t_{n+1}}^{t_{n+1},X^{n+1}}-Y^{n+1}\Big]
    + \frac{\Delta t_n}{2}\left(f_{t_n}^{t_n,X^n}-f^n\right)\\
 & + \frac{\Delta t_n}{2}\mathbb{E}_{t_n}^{X^n}\left[f_{t_{n+1}}^{t_n,X^n}
 - f_{t_{n+1}}^{t_{n+1},X^{n+1}} + f_{t_{n+1}}^{t_{n+1},X^{n+1}}- f^{n+1}\right]
 +R_y^n\\
= \;& \mathbb{E}_{t_n}^{X^n}\left[e^{n+1}_y\right]+\frac{\Delta t_n}{2}
e^n_f+\frac{\Delta t_n}{2} \mathbb{E}_{t_n}^{X^n}[e_f^{n+1}]
+\mathcal{R}_{y_1}^n+\frac{\Delta t_n}{2} \mathcal{R}_{y_2}^n+R^n_y.
\end{aligned}
\end{equation*}
Then under the conditions of the theorem, we have the estimate
\begin{equation*}\label{s5:TL:e7}
\begin{aligned}
\left|e_y^n\right|\leq&\; \left|\mathbb{E}_{t_n}^{X^n}[e^{n+1}_y]\right|+\frac{\Delta t_n}{2} \left|\mathbb{E}_{t_n}^{X^n}[e_f^{n+1}]\right|+\frac{\Delta t_n}{2} |e^n_f|+|\mathcal{R}_{y_1}^n|+\frac{\Delta t_n}{2} |\mathcal{R}_{y_2}^n|+|R^n_y|\\
\leq  &\;\left|\mathbb{E}_{t_n}^{X^n}[e_y^{n+1}]\right|+\frac{\Delta t_n}{2} L \mathbb{E}_{t_n}^{X^n}\left[|e_y^{n+1}|+|e_z^{n+1}|+|e_{\Gamma}^{n+1}|\right]\\
&\;\;+\frac{\Delta t_n }{2} L(|e_y^n|+|e_z^n|+|e_{\Gamma}^n|)+|\mathcal{R}_{y_1}^n|+\frac{\Delta t_n}{2} |\mathcal{R}_{y_2}^n|+|R_y^n|.
\end{aligned}
\end{equation*}
For the squared error $|e_y^n|^2$, given any positive real number $\gamma$
and positive integer $m$, by using the inequalities
$(a+b)^2\leq (1+\gamma\Delta t)a^2
+(1+\frac{1}{\gamma\Delta t})b^2$ and %\;\; \mbox{ and } \;\;
$(\sum_{n=1}^m a_n)^2 \leq m \sum_{n=1}^m a_n^2$,
we have the following estimate.
{\allowdisplaybreaks
\begin{align}
|e_y^n|^2\leq&\; (1+\gamma\Delta t)\left|\mathbb{E}_{t_n}^{X^n}[e_y^{n+1}]\right|^2
+\left(1+\frac{1}{\gamma\Delta t}\right)\Bigg\{\frac{\Delta t_n}{2} L\left(|e_y^n|+|e_z^n|+|e_{\Gamma}^{n}|\right) \nonumber\\
& + \frac{\Delta t_n }{2}L \mathbb{E}_{t_n}^{X^n}\left[|e_y^{n+1}|
+|e_z^{n+1}|+|e_{\Gamma}^{n+1}|\right] +|\mathcal{R}_{y_1}^n|+\frac{\Delta t_n}{2} |\mathcal{R}_{y_2}^n|+|R_y^n|\Bigg\}^2 \nonumber\\
\leq& \; (1+\gamma\Delta t)\left|\mathbb{E}_{t_n}^{X^n}[e_y^{n+1}]\right|^2
+\Bigg\{\frac{15\Delta t^2 L^2}{4}\left(|e_y^n|^2+|e_z^n|^2+|e_{\Gamma}^{n}|^2\right) \label{s5:TL:e8} \\
&+ \frac{15\Delta t^2 L^2}{4}
\mathbb{E}_{t_n}^{X^n}\left[|e_y^{n+1}|^2+|e_z^{n+1}|^2+|e_{\Gamma}^{n+1}|^2\right] +5|\mathcal{R}_{y_1}^n|^2\nonumber\\
& +\frac{5\Delta t^2}{4} |\mathcal{R}_{y_2}^n|^2+5|R_y^n|^2\Bigg\} +\frac{1}{\gamma}\Bigg\{\frac{15\Delta t L^2}{4}(|e_y^n|^2+|e_z^n|^2
+|e_{\Gamma}^{n}|^2) \nonumber\\
& +\frac{15\Delta t L^2 }{4}
\mathbb{E}_{t_n}^{X^n}[|e_y^{n+1}|^2+|e_z^{n+1}|^2+|e_{\Gamma}^{n+1}|^2]\Bigg\} \nonumber\\
& +\frac{1}{\gamma\Delta t}\Bigg\{5|\mathcal{R}_{y_1}^n|^2
+\frac{5\Delta t^2 }{4} |\mathcal{R}_{y_2}^n|^2+5|R_y^n|^2\Bigg\}. \nonumber
\end{align}
}
%\end{equation}

$\bullet$ {\bf  Step 2}:~{\em Estimating the error $e_z^n = {Z}_{t_{n}}^{t_n,X^n}-Z^{n}$.}
%$n=N-1,\dots,1,0$.

Subtracting the scheme (\ref{s4:Z}) from the reference equation (\ref{s3:e10}), we obtain
\begin{equation}\label{s5:TL:e9}
\begin{aligned}
\frac{\Delta t_n}{2}e^n_z
&= \mathbb{E}_{t_n}^{X^n}\left[(Y_{t_{n+1}}^{t_n,X^n}-Y^{n+1})\Delta \tilde{W}_{t_{n+1}}^{\top}\right]\\
& \hspace{0.5cm}+\Delta t_n \mathbb{E}_{t_n}^{X^n}\left[(f_{t_{n+1}}^{t_n,X^n}-f^{n+1})\Delta \tilde{W}_{t_{n+1}}^{\top}\right]+R_z^n.
\end{aligned}
\end{equation}
Substituting the identities
\begin{equation*}
\begin{aligned}
& \mathbb{E}_{t_n}^{X^n}\left[(Y_{t_{n+1}}^{t_n,X^n}-Y^{n+1})\Delta \tilde{W}_{t_{n+1}}^{\top}\right]
= \mathcal{R}_{z_1}^n+\mathbb{E}_{t_n}^{X^n}\left[e_y^{n+1}\Delta \tilde{W}_{t_{n+1}}^{\top}\right],\\
& \mathbb{E}_{t_n}^{X^n}\left[(f_{t_{n+1}}^{t_n,X^n}-f^{n+1})\Delta
\tilde{W}_{t_{n+1}}^{\top}\right]
= \mathcal{R}_{z_2}^n+\mathbb{E}_{t_n}^{X^n}\left[e_f^{n+1}\Delta
\tilde{W}_{t_{n+1}}^{\top}\right],
\end{aligned}
\end{equation*}
into \eqref{s5:TL:e9}, $|e_z^n|$ can be estimated by
\begin{equation}\label{s5:TL:e13}
\begin{aligned}
|e^n_z|=&\left|\frac{2}{\Delta t_n} \mathbb{E}_{t_n}^{X^n}\left[e_y^{n+1}\Delta \tilde{W}_{t_{n+1}}^{\top}\right]
+2\mathbb{E}_{t_n}^{X^n}\left[e_f^{n+1}\Delta \tilde{W}_{t_{n+1}}^{\top}\right]\right.\\
&\hspace{1.5cm}\left. +\frac{2}{\Delta t_n}\mathcal{R}_{z_1}^n+ 2\mathcal{R}_{z_2}^n+\frac{2}{\Delta t_n}R_z^n\right|\\
\leq & \frac{2}{\Delta t_n}\left|\mathbb{E}_{t_n}^{X^n}\left[e_y^{n+1}\Delta \tilde{W}_{t_{n+1}}^{\top}\right]\right|
+ 2\left|\mathbb{E}_{t_n}^{X^n}\left[e_f^{n+1}\Delta \tilde{W}_{t_{n+1}}^{\top}\right]\right|\\
& \hspace{1.5cm}+\frac{2}{\Delta t_n}|\mathcal{R}_{z_1}^n|+ 2|\mathcal{R}_{z_2}^n|+\frac{2}{\Delta t_n}|R_z^n|.
\end{aligned}
\end{equation}
By H\"older's inequality and the inequality
$(a+b)^2\leq (1+\varepsilon)a^2+(1+\frac{1}{\varepsilon})b^2$
for any positive real number $\varepsilon$,
and from \eqref{s5:TL:e13}, we deduce
{\allowdisplaybreaks
\begin{align}
|e^n_z|^2
\leq & (1+\varepsilon)\left(\frac{2}{\Delta t_n}\right)^2
\left|\mathbb{E}_{t_n}^{X^n}\left[e_y^{n+1}\Delta \tilde{W}_{t_{n+1}}^{\top}\right]\right|^2\\
     & + \left(1+\frac{1}{\varepsilon}\right)
       \Bigg\{2\left|\mathbb{E}_{t_n}^{X^n}\left[e_f^{n+1}\Delta \tilde{W}_{t_{n+1}}^{\top}\right]\right|
       + \frac{2|\mathcal{R}_{z_1}^n|}{\Delta t_n}+ 2|\mathcal{R}_{z_2}^n|+\frac{2|R_z^n|}{\Delta t_n}\Bigg\}^2\\
\leq & (1+\varepsilon)\left(\frac{2}{\Delta t_n}\right)^2
\left|\mathbb{E}_{t_n}^{X^n}\left[e_y^{n+1}\Delta \tilde{W}_{t_{n+1}}^{\top}\right]\right|^2
      + 16\left(1+\frac{1}{\varepsilon}\right)
       \Bigg\{\mathbb{E}_{t_n}^{X^n}\left[|e_f^{n+1}|^2\right]
       \mathbb{E}_{t_n}^{X^n}\left[|\Delta \tilde{W}_{t_{n+1}}^{\top}|^2\right]\\
     & + \left(\frac{1}{\Delta t_n}\right)^2|\mathcal{R}_{z_1}^n|^2+|\mathcal{R}_{z_2}^n|^2
       + \left(\frac{1}{\Delta t_n}\right)^2|R_z^n|^2\Bigg\}.\\
\end{align}
}
By the equality $\mathbb{E}_{t_n}^{X^n}[|\Delta \tilde{W}_{t_{n+1}}|^2]=\Delta t_n$, and the estimates of
$\mathbb{E}_{t_n}^{X^n}[|e_f^{n+1}|^2]$ and $|\mathbb{E}_{t_n}^{X^n}[e_y^{n+1}\Delta \tilde{W}_{t_{n+1}}^{\top}]|^2$, i.e.,
\begin{equation}\label{eff}
\begin{aligned}
\mathbb{E}_{t_n}^{X^n}\Big[|e_f^{n+1}|^2\Big]
&\leq \mathbb{E}_{t_n}^{X^n}\left[\left|L(|e_y^{n+1}|+|e_z^{n+1}|+|e_{\Gamma}^{n+1}|)\right|^2\right]\\
&\leq 3L^2\mathbb{E}_{t_n}^{X^n}\Big[|e_y^{n+1}|^2+|e_z^{n+1}|^2+|e_{\Gamma}^{n+1}|^2\Big],
\end{aligned}
\end{equation}
and
\begin{equation}
\begin{aligned}
\left|\mathbb{E}_{t_n}^{X^n}\left[e_y^{n+1}\Delta \tilde{W}_{t_{n+1}}^{\top}\right]\right|^2
=&\; \left|\mathbb{E}_{t_n}^{X^n}\left[(e_y^{n+1}-\mathbb{E}_{t_n}^{X^n}[e_y^{n+1}])\Delta \tilde{W}_{t_{n+1}}^{\top}\right]\right|^2\\
\leq&\; \mathbb{E}_{t_n}^{X^n}\left[|\Delta \tilde{W}_{t_{n+1}}^{\top}|^2\right]
\mathbb{E}_{t_n}^{X^n}\left[(e_y^{n+1}-\mathbb{E}_{t_n}^{X^n}[e_y^{n+1}])^2\right]\\
=&\; \Delta
t_n\left\{\mathbb{E}_{t_n}^{X^n}\left[|e_y^{n+1}|^2\right]-\left|\mathbb{E}_{t_n}^{X^n}[e_y^{n+1}]\right|^2\right\},
\end{aligned}
\end{equation}
into \eqref{s5:TL:e14}, and dividing both sides of the resulting inequality by $(1+\varepsilon)\frac{8}{\Delta t}$,
we obtain an estimate of $|e_z^n|^2$, i.e.,
\begin{equation}\label{s5:TL:e19}
\begin{aligned}
\frac{\Delta t}{8(1+\varepsilon)}|e^n_z|^2
 \leq &\, \frac{c_0}{2}\left\{\mathbb{E}_{t_n}^{X^n}[|e_y^{n+1}|^2]
    - \left|\mathbb{E}_{t_n}^{X^n}[e_y^{n+1}]\right|^2\right\}
    + \frac{6L^2}{\varepsilon}\Delta t^2\mathbb{E}_{t_n}^{X^n}\Big[|e_y^{n+1}|^2\\
& +|e_z^{n+1}|^2+|e_{\Gamma}^{n+1}|^2\Big]
    + \frac{2\Delta t}{\varepsilon}\left\{\frac{|\mathcal{R}_{z_1}^n|^2}{(\Delta t_n)^2}
    + |\mathcal{R}_{z_2}^n|^2+\frac{|R_z^n|^2}{(\Delta t_n)^2}\right\}.
\end{aligned}
\end{equation}

$\bullet$ {\bf  Step 3}:~{\em Estimating the error $e_\Gamma^n = {\Gamma}_{t_{n}}^{t_n,X^n}-\Gamma^{n}$.}

Subtracting the scheme \eqref{s4:G} and the reference equation \eqref{s3:e14}, we have
\begin{equation}\label{s5:TL:e16}
\begin{aligned}
 \frac{1}{2}\Delta t_n e_{\Gamma}^n
=\mathbb{E}_{t_n}^{X^n}\left[(Y_{t_{n+1}}^{t_n,X^n}-Y^{n+1})
\Delta\tilde{\mu}_{t_{n+1}}^*\right]
+\Delta t_n \mathbb{E}_{t_n}^{X^n}\left[(f_{t_{n+1}}^{t_n,X^n}-f^{n+1})
\Delta\tilde{\mu}_{t_{n+1}}^*\right]
+R_{\Gamma}^n.
\end{aligned}
\end{equation}
Substituting
\begin{equation*}
\begin{aligned}
&\mathbb{E}_{t_n}^{X^n}\left[(Y_{t_{n+1}}^{t_n,X^n}-Y^{n+1})
\Delta\tilde{\mu}_{t_{n+1}}^* \right]\\
=&\;\mathbb{E}_{t_n}^{X^n}\left[(Y_{t_{n+1}}^{t_n,X^n}-Y_{t_{n+1}}^{t_{n+1},X^{n+1}})
\Delta\tilde{\mu}_{t_{n+1}}^*\right]+\mathbb{E}_{t_n}^{X^n}\left[e_y^{n+1}\Delta\tilde{\mu}_{t_{n+1}}^*\right]\\
= &\, \mathcal{R}_{\Gamma_1}^n +\mathbb{E}_{t_n}^{X^n}\left[e_y^{n+1}\Delta\tilde{\mu}_{t_{n+1}}^*\right]
\end{aligned}
\end{equation*}
and
\begin{equation*}
\begin{aligned}
& \mathbb{E}_{t_n}^{X^n}\left[(f_{t_{n+1}}^{t_n,X^n}-f^{n+1})
\Delta\tilde{\mu}_{t_{n+1}}^*\right]\\
= &\, \mathbb{E}_{t_n}^{X^n}\left[(f_{t_{n+1}}^{t_n,X^n}-f_{t_{n+1}}^{t_{n+1},X^{n+1}})
\Delta\tilde{\mu}_{t_{n+1}}^*\right]
+\mathbb{E}_{t_n}^{X^n}\left[e_f^{n+1}\Delta\tilde{\mu}_{t_{n+1}}^*\right]\\
= &\, \mathcal{R}_{\Gamma_2}^n +\mathbb{E}_{t_n}^{X^n}\left[e_f^{n+1}
\Delta\tilde{\mu}_{t_{n+1}}^*\right]
\end{aligned}
\end{equation*}
into \eqref{s5:TL:e16}, we obtain an expression of $e_\Gamma^n$ as
\begin{equation}\label{s5:TL:e17}
\begin{aligned}
&e_{\Gamma}^{n}
 = \frac{2}{\Delta t_n}\mathbb{E}_{t_n}^{X^n}\left[e_y^{n+1}\Delta\tilde{\mu}_{t_{n+1}}^*\right]
 + 2\mathbb{E}_{t_n}^{X^n}\left[e_f^{n+1}\Delta\tilde{\mu}_{t_{n+1}}^*\right]\\
&\hspace{0.8cm} + \frac{2}{\Delta t_n}\mathcal{R}_{\Gamma_1}^n + 2\mathcal{R}_{\Gamma_2}^n+\frac{2}{\Delta t_n}R_{\Gamma}^n,
 \end{aligned}
\end{equation}
and consequently we obtain an upper bound of $|e_\Gamma^n|$, i.e.,
\begin{equation}\label{s5:TL:e18}
\begin{aligned}
&|e_{\Gamma}^{n}|
\leq \frac{2}{\Delta t_n}\left|\mathbb{E}_{t_n}^{X^n}[e_y^{n+1}
     \Delta\tilde{\mu}_{t_{n+1}}^*]\right|
 + 2\left|\mathbb{E}_{t_n}^{X^n}[e_f^{n+1}\Delta\tilde{\mu}_{t_{n+1}}^*]\right|\\
 &\hspace{1.5cm}+ \frac{2}{\Delta t_n}|\mathcal{R}_{\Gamma_1}^n|+2|\mathcal{R}_{\Gamma_2}^n|+\frac{2}{\Delta t_n}|R_{\Gamma}^n|.
\end{aligned}
\end{equation}
By H\"older's inequality and the inequality
$(a+b)^2\leq (1+\varepsilon)a^2+(1+\frac{1}{\varepsilon})b^2$
for any positive real number $\varepsilon$,
we obtain the following inequality from \eqref{s5:TL:e18}, i.e.,
\begin{equation}\label{s5:TL:e26}
\begin{aligned}
|e_{\Gamma}^n|^2
\leq &\, (1+\varepsilon)\left(\frac{2}{\Delta t_n}\right)^2
\left|\mathbb{E}_{t_n}^{X^n}\left[e_y^{n+1}\Delta\tilde{\mu}_{t_{n+1}}^*\right]\right|^2\\
     &+ \left(1+\frac{1}{\varepsilon}\right)
       \left\{2\left|\mathbb{E}_{t_n}^{X^n}\left[e_f^{n+1}\Delta\tilde{\mu}_{t_{n+1}}^*\right]\right|
     + \frac{2|\mathcal{R}_{\Gamma_1}^n|}{\Delta t_n}
     + 2|\mathcal{R}_{\Gamma_2}^n|+\frac{2|R_\Gamma^n|}{\Delta t_n}\right\}^2\\
\leq &\, (1+\varepsilon)\left(\frac{2}{\Delta t_n}\right)^2\left|
       \mathbb{E}_{t_n}^{X^n}\left[e_y^{n+1}\Delta\tilde{\mu}_{t_{n+1}}^*\right]\right|^2\\
     & + 16\left(1+\frac{1}{\varepsilon}\right)
        \Bigg\{\mathbb{E}_{t_n}^{X^n}\left[|e_f^{n+1}|^2\right]
           \mathbb{E}_{t_n}^{X^n}\left[|\Delta\tilde{\mu}_{t_{n+1}}^*|^2\right]\\
     & +\left(\frac{1}{\Delta t_n}\right)^2|\mathcal{R}_{\Gamma_1}^n|^2
     + |\mathcal{R}_{\Gamma_2}^n|^2+\Big(\frac{1}{\Delta t_n}\Big)^2|R_{\Gamma}^n|^2\Bigg\}.\\
\end{aligned}
\end{equation}
By substituting the identity
$\mathbb{E}_{t_n}^{X^n}[|\Delta\tilde{\mu}_{t_{n+1}}^*|^2]
=\Delta t_n \int_E \eta^2(e)\lambda(de)$, the estimate of $\mathbb{E}_{t_n}^{X^n}[|e_f^{n+1}|^2]$ given in \eqref{eff},
and the estimate
$|\mathbb{E}_{t_n}^{X^n}[e_y^{n+1}\Delta\tilde{\mu}_{t_{n+1}}^*]|^2$
into \eqref{s5:TL:e26},
and dividing both sides of the resulting inequality by
\(\frac{8(1+\varepsilon)\int_E \rho^2(e)\lambda(de)}{\Delta t}\),  we deduce

\begin{equation}\label{s5:TL:e29}
\begin{aligned}
& \frac{\Delta t}{8(1+\varepsilon)\int_E \eta^2(e)\lambda(de)}|e_{\Gamma}^n|^2 \\
\leq &\; \frac{c_0}{2}\left\{\mathbb{E}_ {t_n}^{X^n}\left[|e_y^{n+1}|^2\right]-\left|\mathbb{E}_{t_n}^{X^n}\left[e_y^{n+1}\right]\right|^2\right\}\\
  &+  \frac{6L^2(\Delta t)^2}{\varepsilon}
       \mathbb{E}_{t_n}^{X^n}\left[|e_y^{n+1}|^2+|e_z^{n+1}|^2+|e_{\Gamma}^{n+1}|^2\right]\\
  & +  \frac{2\Delta t}{\varepsilon \int_E\eta^2(e)\lambda(de)}
       \left[\left(\frac{1}{\Delta t_n}\right)^2|\mathcal{R}_{\Gamma_1}^n|^2
  +   |\mathcal{R}_{\Gamma_2}^n|^2+\left(\frac{1}{\Delta t_n}\right)^2|R_{\Gamma}^n|^2\right].\\
\end{aligned}
\end{equation}

$\bullet$ {\bf  Step 4}:~{\em Combining the estimates from Steps 1-3.}

Now we add \eqref{s5:TL:e8} multiplied by the constant $c_0$,
\eqref{s5:TL:e19}, and \eqref{s5:TL:e29} together to obtain the inequality
{\allowdisplaybreaks
\begin{align}
& c_0|e_y^n|^2 + \frac{\Delta t}{8(1+\varepsilon)}|e^n_z|^2
 + \frac{\Delta t}{8(1+\varepsilon)\int_E \eta^2(e)\lambda(de)}|e_{\Gamma}^n|^2  \label{s5:TL:e88}\\
\leq &\, c_0\left[1+\left(\gamma+\frac{15L^2}{4\gamma}+\frac{15L^2\Delta t}{4}
 + \frac{12L^2\Delta t}{c_0\varepsilon}\right)\Delta t\right]\mathbb{E}_{t_n}^{X^n}\left[|e_y^{n+1}|^2\right]\nonumber\\
& + \left[\frac{15c_0}{4\gamma}+\left(\frac{15c_0}{4}+\frac{12}{\varepsilon}\right)
 \Delta t\right]L^2\Delta t\mathbb{E}_{t_n}^{X^n}\left[|e_z^{n+1}|^2\right]\nonumber\\
& + \left[\frac{15c_0}{4\gamma}+\left(\frac{15c_0}{4}+\frac{12}{\varepsilon}\right)
 \Delta t\right]L^2\Delta t\mathbb{E}_{t_n}^{X^n}\left[|e_{\Gamma}^{n+1}|^2\right]\nonumber\\
&+\left(\frac{15c_0}{4\gamma}+\frac{15c_0\Delta t}{4}\right)
 L^2\Delta t \left(|e_y^n|^2+|e_z^n|^2+|e_{\Gamma}^n|^2\right)\nonumber\\
& + 5c_0\left(1+\frac{1}{\gamma\Delta t}\right) \left\{|\mathcal{R}_{y_1}^n|^2+\frac{1}{4}\Delta t^2 |\mathcal{R}_{y_2}^n|^2+|R_y^n|^2\right\}\nonumber\\
& +\frac{2\Delta t}{\varepsilon}
\left\{\frac{1}{(\Delta t_n)^2}|\mathcal{R}_{z_1}^n|^2+|\mathcal{R}_{z_2}^n|^2+\frac{1}{(\Delta t_n)^2}|R_z^n|^2\right\}\nonumber\\
& +\frac{2\Delta t}{\varepsilon\int_{E}\eta^2(e) \lambda(de)}
\left\{\frac{1}{(\Delta t_n)^2}|\mathcal{R}_{\Gamma_1}^n|^2+|\mathcal{R}_{\Gamma_2}^n|^2+\frac{1}{(\Delta t_n)^2}|R_{\Gamma}^n|^2\right\}.\nonumber
\end{align}}
By taking expectation $\mathbb{E}[\cdot]$ on both sides of \eqref{s5:TL:e88}, we deduce
%\begin{equation}
\begin{align}
& c_0\left(1-C_1\Delta t\right)\mathbb{E}\left[|e_y^n|^2\right]+C_3\Delta t \mathbb{E}\big[|e^n_z|^2\big]
+C_6\Delta t\mathbb{E}\big[|e_{\Gamma}^n|^2\big]\label{s5:TL:e25}\\
 \leq &\, c_0(1+C_2\Delta t) \mathbb{E}\big[|e_y^{n+1}|^2\big]
      +  C_4\Delta t \mathbb{E}\big[|e^{n+1}_z|^2\big]
      +  C_4\Delta t \mathbb{E}\big[|e^{n+1}_{\Gamma}|^2\big]\nonumber\\
   & + \frac{C_5}{\Delta t}\mathbb{E}\left[|\mathcal{R}_{y_1}^n|^2
               +\frac{1}{4}\Delta t^2 |\mathcal{R}_{y_2}^n|^2+|R_y^n|^2\right] \nonumber\\
& +\frac{2\Delta t}{\varepsilon}\mathbb{E}\left[\left(\frac{1}{\Delta t_n}\right)^2|\mathcal{R}_{z_1}^n|^2
               +|\mathcal{R}_{z_2}^n|^2+\left(\frac{1}{\Delta t_n}\right)^2|R_{z}^n|^2\right]\nonumber\\
& + \frac{2\Delta t}{\varepsilon\int_{E}\eta^2(e) \lambda(de)}
                \mathbb{E}\left[\frac{1}{(\Delta t_n)^2}|\mathcal{R}_{\Gamma_1}^n|^2+|\mathcal{R}_{\Gamma_2}^n|^2
               +\frac{1}{(\Delta t_n)^2}|R_{\Gamma}^n|^2\right],\nonumber
\end{align}
%\end{equation}
%
where the constants $C_1$,  $C_2$, $C_3$, $C_4$, $C_5$, $C_6$ are defined by
\begin{equation}
\begin{aligned}
&C_1=\left(\frac{15}{4\gamma}+\frac{15\Delta t}{4}\right)L^2,\; C_2=\left(\gamma+\frac{15L^2}{4\gamma}+\frac{15L^2\Delta
t}{4}+\frac{12L^2\Delta t}{c_0\varepsilon}\right),\\
&C_3=\frac{1}{8(1+\varepsilon)}-\left(\frac{15c_0}{4\gamma}+\frac{15c_0\Delta t}{4}\right)L^2,\; C_4=\left[\frac{15c_0}{4\gamma}+\left(\frac{15c_0}{4}+\frac{12}{\varepsilon}\right)\Delta t\right]L^2,\\
& C_5=5c_0\frac{1+\gamma\Delta t}{\gamma}, \; C_6=\frac{1}{8(1+\varepsilon)\int_E\eta^2(e)\lambda(de)}-\left(\frac{15c_0}{4\gamma}+\frac{15c_0\Delta t}{4}\right)L^2.
\end{aligned}
\end{equation}

Now we set $\varepsilon=1$, $\gamma$ large enough and $\Delta
t$ sufficiently small, such that if $0<\Delta t\leq \Delta t$
then $C_1\le C$, $C_2\le C$, $C_5\le C$, $1-C\Delta t>0$, and
$C_3-C_4>C^*>0$, $C_6-C_4>C^*>0$ where $C$ and $C^*$ are two positive constants
depending on $c_0$ and $L$. Then for $0<\Delta t\leq \Delta t_0$,
 we deduce from \eqref{s5:TL:e25}
{\allowdisplaybreaks
\begin{align*}
& c_0(1-C\Delta t)\mathbb{E}\big[|e_y^n|^2\big]+C_3\Delta t\mathbb{E}\big[|e^n_z|^2\big]
 + C_6\Delta t\mathbb{E}\big[|e_{\Gamma}^n|^2\big]\\
 %%
% \leq &\, c_0(1+C\Delta t)\mathbb{E}\big[|e_y^{n+1}|^2\big]+C_4\Delta t \mathbb{E}\big[|e^{n+1}_z|^2\big]
% + C_4\Delta t \mathbb{E}\big[|e^{n+1}_{\Gamma}|^2\big]\\
% %
%& + \frac{C}{\Delta t}{\mathbb{E}\left[|\mathcal{R}_{y_1}^n|^2+(\Delta t)^2 |\mathcal{R}_{y_2}^n|^2+|R_y^n|^2\right]}\\
%%
%& + C\Delta t\mathbb{E}\left[\left(\frac{1}{\Delta t_n}\right)^2|\mathcal{R}_{z_1}^n|^2
% + |\mathcal{R}_{z_2}^n|^2+\left(\frac{1}{\Delta t_n}\right)^2|R_{z}^n|^2\right]\\
% %
%& +C \Delta t\mathbb{E}\left[\left(\frac{1}{\Delta t_n}\right)^2|\mathcal{R}_{\Gamma_1}^n|^2
% + |\mathcal{R}_{\Gamma_2}^n|^2+\left(\frac{1}{\Delta t_n}\right)^2|R_{\Gamma}^n|^2\right].\\
% %
 \leq &\, c_0(1+C\Delta t)\mathbb{E}\big[|e_y^{n+1}|^2\big]+C_4\Delta t \mathbb{E}\big[|e^{n+1}_z|^2\big]
 + C_4\Delta t \mathbb{E}\big[|e^{n+1}_{\Gamma}|^2\big]\\
& + C\mathbb{E}\left[\frac{1}{\Delta t}|\mathcal{R}_{y_1}^n|^2
 + \Delta t |\mathcal{R}_{y_2}^n|^2+\frac{1}{\Delta t}|R_y^n|^2\right]\\
& + C\mathbb{E}\left[\frac{1}{\Delta t}|\mathcal{R}_{z_1}^n|^2
 + \Delta t|\mathcal{R}_{z_2}^n|^2+\frac{1}{\Delta t}|R_{z}^n|^2\right]\\
& +C \mathbb{E}\left[\frac{1}{\Delta t}|\mathcal{R}_{\Gamma_1}^n|^2
 + \Delta t|\mathcal{R}_{\Gamma_2}^n|^2+\frac{1}{\Delta t}|R_{\Gamma}^n|^2\right].
\end{align*}}
Dividing both sides of the upper inequality by $(1-C\Delta t)$, we easily get
\begin{equation}\label{s5:TL:e22}
\begin{aligned}
& c_0\mathbb{E}[|e_y^n|^2]+C_3\Delta t
\mathbb{E}[|e^n_z|^2]+C_6\Delta t\mathbb{E}[|e_{\Gamma}^n|^2]\\
 \leq &\, \frac{1+C\Delta t}{1-C\Delta t}\Big(c_0\mathbb{E}[|e_y^{n+1}|^2]
          +   C_4\Delta t \mathbb{E}[|e^{n+1}_z|^2]
          +   C_4\Delta t \mathbb{E}[|e^{n+1}_{\Gamma}|^2\Big)\\
& + \frac{C}{1-C\Delta t}\mathbb{E}\left[\frac{1}{\Delta t}|\mathcal{R}_{y_1}^n|^2
          +\Delta t |\mathcal{R}_{y_2}^n|^2+\frac{1}{\Delta t}|R_y^n|^2\right]\\
& + \frac{C}{1-C\Delta t}\mathbb{E}\left[\frac{1}{\Delta t}|\mathcal{R}_{z_1}^n|^2
          + \Delta t|\mathcal{R}_{z_2}^n|^2+\frac{1}{\Delta t}|R_{z}^n|^2\right]\\
& + \frac{C}{1-C\Delta t}\mathbb{E}\left[\frac{1}{\Delta t}|\mathcal{R}_{\Gamma_1}^n|^2
          + \Delta t|\mathcal{R}_{\Gamma_2}^n|^2+\frac{1}{\Delta t}|R_{\Gamma}^n|^2\right].\\
\end{aligned}
\end{equation}
From the inequality \eqref{s5:TL:e22}, by recursively inserting
$e_y^{i}$, $i=n+1,\dots,N-1$, we deduce
{\allowdisplaybreaks
\begin{align}
& c_0\mathbb{E}\left[|e_y^n|^2\right]+C_3\Delta t\sum\limits_{i=n}^{N-1}
 \left(\frac{1+C\Delta t}{1-C\Delta t}\right)^{i-n}\mathbb{E}\left[|e_z^i|^2\right] \nonumber
+ C_6\Delta t\sum\limits_{i=n}^{N-1}\left(\frac{1+C\Delta t}{1-C\Delta t}\right)^{i-n}
   \mathbb{E}\left[|e_{\Gamma}^i|^2\right] \nonumber\\
\leq &\, \left(\frac{1+C\Delta t}{1-C\Delta t}\right)^{N-n}c_0\mathbb{E}\left[|e_y^N|^2\right]
 + C_4\Delta t\sum\limits_{i=n+1}^N \left(\frac{1+C\Delta t}{1-C\Delta t}\right)^{i-n}\mathbb{E}\left[|e_z^i|^2\right]\nonumber\\
& + C_4\Delta t\sum\limits_{i=n+1}^N
   \left(\frac{1+C\Delta t}{1-C\Delta t}\right)^{i-n}\mathbb{E}[|e_{\Gamma}^i|^2]\nonumber\\
& + \sum\limits_{i=n}^{N-1}\left(\frac{1+C\Delta t}{1-C\Delta t}\right)^{i-n}
   \frac{C}{1-C\Delta t}\mathbb{E}\left[\frac{1}{\Delta t}|\mathcal{R}_{y_1}^i|^2
   + \Delta t |\mathcal{R}_{y_2}^i|^2+\frac{1}{\Delta t}|R_y^i|^2\right]\nonumber\\
& + \sum\limits_{i=n}^{N-1}\left(\frac{1+C\Delta t}{1-C\Delta t}\right)^{i-n}
   \frac{C}{1-C\Delta t}\mathbb{E}\left[\frac{1}{\Delta t}|\mathcal{R}_{z_1}^i|^2
   + \Delta t|\mathcal{R}_{z_2}^i|^2+\frac{1}{\Delta t}|R_{z}^i|^2\right]\nonumber\\
& + \sum\limits_{i=n}^{N-1}\left(\frac{1+C\Delta t}{1-C\Delta t}\right)^{i-n}
   \frac{C}{1-C\Delta t}\mathbb{E}\left[\frac{1}{\Delta t}|\mathcal{R}_{\Gamma_1}^i|^2
   + \Delta t|\mathcal{R}_{\Gamma_2}^i|^2+\frac{1}{\Delta t}|R_{\Gamma}^i|^2\right], \nonumber
\end{align}}
which immediately leads to
{\allowdisplaybreaks
\begin{align}
& c_0\mathbb{E}\left[|e_y^n|^2\right]+ C^* \Delta t\sum\limits_{i=n}^{N-1}
 \left(\frac{1+C\Delta t}{1-C\Delta t}\right)^{i-n}\mathbb{E}\left[|e_z^i|^2\right]
 + C^*\Delta t\sum\limits_{i=n}^{N-1}\left(\frac{1+C\Delta t}{1-C\Delta t}\right)^{i-n}\mathbb{E}\left[|e_{\Gamma}^i|^2\right]\nonumber\\
\leq &\, \left(\frac{1+C\Delta t}{1-C\Delta t}\right)^{N-n}
   c_0\mathbb{E}\left[|e_y^N|^2\right]+C_4\Delta t\left(\frac{1+C\Delta t}{1-C\Delta t}\right)^{N-n} \mathbb{E}\left[|e_z^N|^2\right]\nonumber\\
& + C_4\Delta t\left(\frac{1+C\Delta t}{1-C\Delta t}\right)^{N-n} \mathbb{E}\left[|e_{\Gamma}^N|^2\right]\nonumber\\
& + \sum\limits_{i=n}^{N-1}\left(\frac{1+C\Delta t}{1-C\Delta t}\right)^{i-n}
  \frac{C}{1-C\Delta t}\mathbb{E}\left[\frac{1}{\Delta t}|\mathcal{R}_{y_1}^i|^2
  + \Delta t |\mathcal{R}_{y_2}^i|^2+\frac{1}{\Delta t}|R_y^i|^2\right]\nonumber\\
& + \sum\limits_{i=n}^{N-1}\left(\frac{1+C\Delta t}{1-C\Delta t}\right)^{i-n}
    \frac{C}{1-C\Delta t}\mathbb{E}\left[\frac{1}{\Delta t}|\mathcal{R}_{z_1}^i|^2
  + \Delta t|\mathcal{R}_{z_2}^i|^2+\frac{1}{\Delta t}|R_{z}^i|^2\right]\nonumber\\
& + \sum\limits_{i=n}^{N-1}\left(\frac{1+C\Delta t}{1-C\Delta t}\right)^{i-n}
    \frac{C}{1-C\Delta t}\mathbb{E}\left[\frac{1}{\Delta t}|\mathcal{R}_{\Gamma_1}^i|^2
  + \Delta t|\mathcal{R}_{\Gamma_2}^i|^2+\frac{1}{\Delta t}|R_{\Gamma}^i|^2\right].\nonumber
\end{align}}
The proof is completed. \end{proof}

\begin{rem}
It is worth to note that Theorem \ref{s5:TL} also implies that Scheme \ref{s4:sch:1} is stable with respect to the terminal condition, that is,
for any $\varepsilon >0$, there exists a positive number $\delta >0$, such that, if
 $\mathbb{E}[|\overline{Y}^N-Y^N|^2]<\delta$,
$\mathbb{E}[|\overline{Z}^N-Z^N|^2]<\delta$ and
$\mathbb{E}[|\overline{\Gamma}^N-\Gamma^N|^2]<\delta$, where
$(\overline{Y}^N,\overline{Z}^N,\overline{\Gamma}^N)$ and $(Y^N,Z^N,\Gamma^N)$
are two different terminal conditions,
then for $0\leq n\leq N-1$, it holds that
$$\mathbb{E}[|\overline{Y}^n-Y^n|^2]
+\Delta t \sum\limits_{i=n}^{N}\mathbb{E}[|\overline{Z}^i-Z^i|^2]
+\Delta t \sum\limits_{i=n}^{N}\mathbb{E}[|\overline{\Gamma}^i-\Gamma^i|^2]<\varepsilon.$$
\end{rem}

The next task is to estimate all the residual terms in \eqref{s5:TL:e1},
and the main technique used is the It\^o-Taylor expansion \cite{Platen:2010eo}.
Under some reasonable regularity conditions on
the data $b$, $\sigma$, $c$, $f$ and $\varphi$ in the FBSDE,
we now derive estimates for the local truncation errors $R_y^n$, $R_{z}^n$ and $R_{\Gamma}^n$ defined
in \eqref{s3:e6}, \eqref{s3:e10} and \eqref{s3:e14}, respectively.
To proceed, we need the following standard assumption.
\begin{assumption}\label{s5:hy}
Under the condition that $X_{0}$ is $\mathcal{A}_{0}$-measurable as well as
$\mathbb{E}[|X_{0}|^2]<\infty$, we assume that $b$, $\sigma$ and
$c$ are jointly $L^2$-measurable in $(t,x)\in[0,T]\times
\mathbb{R}^q$, and there exist real constants $L>0$ and $K>0$
such that
\begin{equation}\label{lip1}
\begin{aligned}
\hspace{0.5cm} &|b(t,x)-b(t,x')|\leq L|x-x'|,\;\;  |\sigma(t,x)-\sigma(t,x')|\leq L|x-x'|,\\
& \int_E |c(t,x,e)-c(t,x',e)|^2 \lambda (de)\leq L|x-x'|^2,
\end{aligned}
\end{equation}
and
\begin{equation}\label{lin1}
\hspace{-0.47cm}\begin{aligned}
& |b(t,x)|^2\leq K(1+|x|^2),\;\; |\sigma(t,x)|^2\leq K(1+|x|^2),\\
 &\int_E |c(t,x,e)|^2 \lambda(de)\leq K(1+|x|^2),\\
\end{aligned}
\end{equation}
 for all $t \in [0,T]$ and $x,x'\in {\mathbb R}^q$.
\end{assumption}
Under Assumption \ref{s5:hy}, if $\mathbb{E}[|X_{0}|^{2m}]<\infty$ for
some integer $m\geq 1$, the solution of the forward SDE in (\ref{s1:e1}) has the
estimate
\begin{equation}\label{CXn}
\mathbb{E}_{t_n}^{X^n}\left[|X_s^{t_n,X^n}|^{2m}\right]\leq
\left(1+\mathbb{E}_{t_n}^{X^n}\big[|X^n|^{2m}\big]\right)\mathrm{e}^{C(s-t_n)},
\end{equation}
where $s\in[t_n,T]$ and $C$ is a positive constant
depending only on the constants $K$, $L$ and $m$.

For the sake of presentation simplicity, in the following lemmas and theorems,
we only consider the one-dimensional case $(q=d=1)$, but our results can be extended to
multidimensional cases without any essential difficulty.
To proceed, we define three partial integro-differential operators:
\begin{equation}\label{LLL}
\begin{aligned}
L^0 v(t,x):= &\, \frac{\partial v}{\partial t}(t,x)+b(t,x)\frac{\partial v}{\partial x}(t,x)
+\frac{1}{2}\sigma^2(t,x) \frac{\partial^2 v}{\partial x^2}(t,x)\\
& +\int_E \left[v(t,x+c(t,x,e))-v(t,x)-\frac{\partial v}{\partial x}(t,x)c(t,x,e)\right]\lambda (de),\\
L^1v(t,x) :=&\, \sigma(t,x)\frac{\partial v}{\partial x}(t,x),\\
L^{-1}v(t,x):=&\, v(t,x+c(t,x,e))-v(t,x);
\end{aligned}
\end{equation}
and introduce the following notation:
\[
\begin{aligned}
& \mathcal{C}_b^{(k_1, \ldots, k_J)}(D_1 \times \cdots \times D_J)\\
:= & \Bigg\{\phi: \prod_{j=1}^J D_j \rightarrow \mathbb{R}\; \bigg|\; \frac{\partial^{\alpha_1}\cdots \partial^{\alpha_J} \phi}{\partial^{\alpha_1} x_1\cdots \partial^{\alpha_J} x_J}\; \text{is bounded and continuous }\\
&\hspace{0.8cm} \text{ for } \alpha_j \le k_j, j = 1, \ldots, J, \text{ where }  \vec{\alpha} := (\alpha_1, \ldots, \alpha_J) \in \mathbb{N}^J \Bigg\},
\end{aligned}
\]
where $J \in \mathbb{N}^{+}$, $D_1 \times \cdots \times D_J \subset \mathbb{R}^J$.

Now we give the estimates of $R_y^n$, $R_{z,1}^n$, $R_{z,2}^n$, $R_{\Gamma,1}^n$ and $R_{\Gamma,2}^n$
in the following Lemma \ref{s5:lem1}.

\begin{lem}\label{s5:lem1}
Under Assumption \ref{s5:hy}, if the data of the FBSDEs in \eqref{s1:e1} satisfy the following regularity conditions: $f(t,x,y,z,\gamma)\in \mathcal{C}_b^{(2,4,4,4,4)}([0,T] \times \mathbb{R}^4)$,
$b(t,x)\in \mathcal{C}_b^{(2,4)}([0,T]\times \mathbb{R})$,
$\sigma(t,x) \in\mathcal{C}_b^{(2,4)}([0,T]\times \mathbb{R})$, $\varphi(x)\in \mathcal{C}_b^{6+\alpha}(\mathbb{R})$
with $\alpha\in(0,1)$, and $c(t,x,e) \in \mathcal{C}_b^{(2,4,\infty)}([0,T]\times \mathbb{R}^2)$, then
for sufficiently small $\Delta t = \max_{n} \Delta t_n$, we have the estimates
\[
\mathbb{E}\left[|R_y^n|^2\right]
\leq C\left(1+\mathbb{E}\left[|X^n|^8\right]\right)(\Delta t)^6,
\]
$$ \mathbb{E}\left[|R_{z,1}^n|^2\right]
 \leq C(1+\mathbb{E}\left[|X^n|^8\right])(\Delta t_n)^6,
\;\;\; \mathbb{E}\left[|R_{z,2}^n|^2\right]
\leq  C(1+\mathbb{E}\left[|X^n|^8\right])(\Delta t_n)^6,
$$
$$
 \mathbb{E}\left[|R_{\Gamma,1}^n|^2\right]
 \leq C(1+\mathbb{E}\left[|X^n|^8\right])(\Delta t_n)^6,
\;\;\; \mathbb{E}\left[|R_{\Gamma,2}^n|^2\right]
\leq  C(1+\mathbb{E}\left[|X^n|^8\right])(\Delta t_n)^6,
$$
where $R_y^n$, $R_{z,1}^n$, $R_{z,2}^n$, $R_{\Gamma,1}^n$ and $R_{\Gamma,2}^n$
are defined in \eqref{s3:Ry}, \eqref{s3:e8}, \eqref{s3:e9}, \eqref{s3:e12} and \eqref{s3:e13}, respectively,
and $C$ is a positive constant depending only on $T$, $K$ and
upper bounds of the derivatives of $b$, $\sigma$, $c$, $f$ and $\varphi$.
\end{lem}

\begin{proof}
Based on the relation shown in \eqref{FK} between the solution $(Y_t, Z_t, \Gamma_t)$
of the BSDE in \eqref{s3:e1} and the solution $u(t,x)$ of the PIDE in \eqref{eq:Kol},
it is easy to prove that $u(t,x) \in \mathcal{C}^{(2,4)}_b([0,T]\times \mathbb{R}^q)$ under the regularity conditions given in \cite{Anonymous:fk}  on $f, b, \sigma, c$ and $\varphi$. Then, for $t \le s \le T$, the function $F=F(t,x)$ defined by
\begin{equation}\label{FF}
F(s,x):=f\left(s,x,u(s,x),\nabla u(s,x)\sigma(s,x), \Gamma(s,x)\right),
\end{equation}
 is in the space $\mathcal{C}_b^{(2,4)}([0,T]\times \mathbb{R}^q)$.
Setting $x = X_{s}^{t_{n},X^{n}}$ in \eqref{FF} and applying the
It\^{o}-Taylor expansion to $F(s,X_s^{t_n,X^n})$, we obtain
\begin{equation}\label{s5:Lem1:e0}
\begin{aligned}
F(s,X_s^{t_n,X^n})
= &\, F(t_n,X^n)+\int_{t_n}^s L^0F(r,X_r^{t_n,X^n})dr
    + \int_{t_n}^s L^1F(r,X_r^{t_n,X^n})dW_r\\
  & +\int_{t_n}^s \int_E L^{-1}F(r,X_{r-}^{t_n,X^n})\tilde{\mu}(de,dr),
\end{aligned}
\end{equation}
where the operators $L^0$, $L^1$ and $L^{-1}$ are defined in \eqref{LLL}.

Taking the conditional expectation $\mathbb{E}_{t_n}^{X^n}[\cdot]$ on both sides,
we have that
$$\mathbb{E}_{t_n}^{X^n}\Big[\int_{t_n}^s L^1F(r,X_r^{t_n,X^n})dW_r\Big]=0\quad \text{and}\quad  \mathbb{E}_{t_n}^{X^n}\Big[\int_{t_n}^s \int_E L^{-1}F(r,X_{r-}^{t_n,X^n})\tilde{\mu}(de,dr)\Big]=0.$$
Then following the same arguments used in the proof of Lemmas 4.2-4.4 in \cite{Zhao:2014eu},
we obtain the estimates of the lemma. The proof is completed.
\end{proof}

From Lemma \ref{s5:lem1} and the definitions of $R_z^n$ in \eqref{s3:e10} and $R_{\Gamma}^n$
in \eqref{s3:e14}, we easily get the estimates of $R_z^n$ and $R_{\Gamma}^n$, stated
in the following lemma.

\begin{lem}\label{s5:lem4}
Under the conditions of Lemma \ref{s5:lem1}, for
sufficiently small time step size $\Delta t = \max_{n} \Delta t_n$, we have that
$$\mathbb{E}\left[|R_z^n|^2\right] \leq C\left(1+\mathbb{E}\left[|X^n|^8\right]\right)(\Delta t)^6\;\text{ for } \;0\leq n\leq N-1,$$
$$\mathbb{E}\left[|R_\Gamma^n|^2\right] \leq C\left(1+\mathbb{E}\left[|X^n|^8\right]\right)(\Delta t)^6\;\text{ for } \;0\leq n\leq N-1,$$
where $R_z^n$, $R_{\Gamma}^{n}$ are defined in \eqref{s3:e10}, \eqref{s3:e14},
$C$ is a positive constant depending on $T$, $K$ and the
upper bounds of the derivatives of $b$, $\sigma$, $c$, $f$ and $\varphi$.
\end{lem}
%\begin{proof} From the definitions of $R_z^n$ in \eqref{s3:e10} and $R_{\Gamma}^n$
%in \eqref{s3:e14}, the lemma is the direct consequence of
%Lemma \ref{s5:lem1}. It can be proved by following the same argument as $R_z^n$, so
%we omit the proof of $R_{\Gamma}^{n}$ here.
%\end{proof}
%

%\begin{lem}\label{s5:lem8}
%Under the conditions of Lemma \ref{s5:lem1}, for
%sufficiently small time step size $\Delta t = \max_{n} \Delta t_n$, we have that
%$$\mathbb{E}\left[|R_\Gamma^n|^2\right] \leq C\left(1+\mathbb{E}\left[|X^n|^8\right]\right)(\Delta t)^6\;\text{ for } \;0\leq n\leq N-1,$$ where $R_{\Gamma}^{n}$ is defined in \eqref{s3:e14},
%$C$ is a positive constant depending on $T$, $K$ and the
%upper bounds of the derivatives of $b$, $\sigma$, $f$ and $\varphi$.
%\end{lem}
%
%\begin{proof}
%It can be proved by following the same argument as in Lemma \ref{s5:lem4}, so
%we omit the proof here.
%\end{proof}

Now combining Theorem \ref{s5:TL}, Lemma \ref{s5:lem1} and Lemma \ref{s5:lem4}
as well as the estimates given in \eqref{s5:st} and \eqref{s5:EX},
we obtain the convergence rate of Scheme \ref{s4:sch:1} in the following theorem.

\begin{theorem}\label{s5:TL2}
Under Lemma \ref{s5:lem1} and Lemma \ref{s5:lem4}, if \eqref{s5:st} and \eqref{s5:EX}
hold for the scheme \eqref{s4:X111} for the forward SDE,
then, for sufficiently small time step size $\Delta t = \max_{n} \Delta t_n$,
the errors $e_y^n$, $e_z^n$ and $e_\Gamma^n$ in \eqref{err} for $n = 0, \ldots, N-1$ can be bounded by
\begin{equation*}
\begin{aligned}
  & \mathbb{E}[|e_y^n|^2]+\Delta t\sum\limits_{i=n}^{N-1}
\left(\frac{1+C\Delta t}{1-C\Delta t}\right)^{i-n} \mathbb{E}[|e^i_z|^2+|e^i_\Gamma|^2]\\
 \leq & C_1(\mathbb{E}[|e_y^{N}|^2]+\Delta t \mathbb{E}[|e^{N}_z|^2+|e^{N}_\Gamma|^2])
 +C_2\Big((\Delta t)^{2\alpha}+(\Delta t)^{2\beta}+(\Delta t)^{2\gamma}+(\Delta t)^4\Big),
\end{aligned}
\end{equation*}
where $\alpha$, $\beta$, $\gamma$ are defined in \eqref{s5:EX},
$C >0$ depends on $c_0$ and $L$,
$C_1>0$ depends on $c_0$, $T$ and $L$, $C_2>0$
depends on $c_0$, $T$, $L$, $K$, $X_0$ and the upper bounds of
the derivatives of $b$, $\sigma$, $c$, $f$ and $\varphi$.
\end{theorem}
\begin{proof}
From the definitions of $\mathcal{R}_{y_1}^i$, $\mathcal{R}_{y_2}^i$,
$\mathcal{R}_{z_1}^i$, $\mathcal{R}_{z_2}^i$, $\mathcal{R}_{\Gamma_1}^i$ and
$\mathcal{R}_{\Gamma_2}^i$ in Theorem \ref{s5:TL},
under the conditions of the theorem, we easily get the estimates
\begin{equation}\label{s5:CL:e1}
\begin{aligned}
\mathbb{E}[|X^i|^2]\leq  &\; C(1+\mathbb{E}[|X_0|^2]),\\
\mathbb{E}[|\mathcal{R}_{y_1}^i|^2]
\leq &\; C(1+\mathbb{E}[|X^i|^{4r_1}])(\Delta t)^{2\beta+2}\le C(1+\mathbb{E}[|X_0|^{4r_1}])(\Delta t)^{2\beta+2},\\
\mathbb{E}[|\mathcal{R}_{y_2}^i|^2]\leq & C(1+\mathbb{E}[|X^i|^{4r_1}])(\Delta t)^{2\beta+2}
\leq\; C(1+\mathbb{E}[|X_0|^{4r_1}])(\Delta t)^{2\beta+2},\\
\mathbb{E}[|\mathcal{R}_{z_1}^i|^2] \leq & C(1+\mathbb{E}[|X^i|^{4r_2}])(\Delta t)^{2\gamma+2}
\leq\; C(1+\mathbb{E}[|X_0|^{4r_2})(\Delta t)^{2\gamma+2},\\
\mathbb{E}[|\mathcal{R}_{z_2}^i|^2] \leq &
C(1+\mathbb{E}[|X^i|^{4r_2}])(\Delta t)^{2\gamma+2} \leq\;
C(1+\mathbb{E}[|X_0|^{4r_2}])(\Delta t)^{2\gamma+2},\\
\mathbb{E}[|\mathcal{R}_{\Gamma_1}^i|^2] \leq & C(1+\mathbb{E}[|X^i|^{4r_3})(\Delta t)^{2\alpha+2}
\leq\; C(1+\mathbb{E}[|X_0|^{4r_3})(\Delta t)^{2\alpha+2},\\
\mathbb{E}[|\mathcal{R}_{\Gamma_2}^i|^2] \leq &
C(1+\mathbb{E}[|X^i|^{4r_3}])(\Delta t)^{2\alpha+2} \leq\;
C(1+\mathbb{E}[|X_0|^{4r_3}])(\Delta t)^{2\alpha+2}
\end{aligned}
\end{equation}
for $i=0,1,\dots,N-1$.
By Lemmas \ref{s5:lem1}, \ref{s5:lem4} and \ref{s5:lem4}, and inequality \eqref{s5:st},
for $0\leq i\leq N-1$, we have
\begin{equation}\label{s5:CL:e2}
\begin{aligned}
\mathbb{E}[|R_y^i|^2]
& \leq  C(1+\mathbb{E}[|X_0|^8])(\Delta t)^6, \quad
\mathbb{E}[|R_z^i|^2]
\leq C(1+\mathbb{E}[|X_0|^8])(\Delta t)^6,\\
%\text{and} \quad
\mathbb{E}[|R_{\Gamma}^i|^2]
&\leq C(1+\mathbb{E}[|X_0|^8])(\Delta t)^6.
\end{aligned}
\end{equation}
By (\ref{s5:CL:e1}) and (\ref{s5:CL:e2}), we deduce
\begin{equation}\label{s5:CL:e3}
\begin{aligned}
&\sum\limits_{i=n}^{N-1}\Big(\frac{1+C\Delta t}{1-C\Delta
t}\Big)^{i-n}\frac{C\mathbb{E}[|\mathcal{R}_{y_1}^i|^2+(\Delta t)^2
|\mathcal{R}_{y_2}^i|^2+|R_y^i|^2]}{\Delta t(1-C\Delta t)}\\
\leq &
C(1+\mathbb{E}[|X_0|^{4r_1}]+\mathbb{E}[|X_0|^8])((\Delta
t)^{2\beta}+(\Delta t)^4),
\end{aligned}
\end{equation}
\begin{equation}\label{s5:CL:e4}
\begin{aligned}
&\sum\limits_{i=n}^{N-1}\Big(\frac{1+C\Delta t}{1-C\Delta t}\Big)^{i-n}
\frac{C\mathbb{E}[|\mathcal{R}_{z_1}^i|^2+(\Delta t)^2
|\mathcal{R}_{z_2}^i|^2+|R_z^i|^2]}{\Delta t(1-C\Delta t)}\\
\leq &
C(1+\mathbb{E}[|X_0|^{4r_2}]+\mathbb{E}[|X_0|^8])((\Delta
t)^{2\gamma}+(\Delta t)^4),
\end{aligned}
\end{equation}
and
\begin{equation}\label{s5:CL:e5}
\begin{aligned}
&\sum\limits_{i=n}^{N-1}\Big(\frac{1+C\Delta t}{1-C\Delta t}\Big)^{i-n}
\frac{C\mathbb{E}[|\mathcal{R}_{\Gamma_1}^i|^2+(\Delta t)^2
|\mathcal{R}_{\Gamma_2}^i|^2+|R_\Gamma^i|^2]}{\Delta t(1-C\Delta t)}\\
\leq &
C(1+\mathbb{E}[|X_0|^{4r_3}]+\mathbb{E}[|X_0|^8])((\Delta
t)^{2\alpha}+(\Delta t)^4).
\end{aligned}
\end{equation}
By applying \eqref{s5:CL:e3}, \eqref{s5:CL:e4}
and \eqref{s5:CL:e5} to Theorem \ref{s5:TL}, we complete the proof.
\end{proof}

\section{The fully discrete scheme}\label{sec:fully}
In this section, we will develop a fully discrete scheme based on Scheme \ref{s4:sch:1}
by assuming that the jump process of $X_t$ in \eqref{s1:e1} has {\em finite activity}.
This means the Poisson random measure $\tilde{\mu}(de, dt)$ can be represented by
\begin{equation}\label{comp}
\tilde{\mu}(de, dt) = \mu(de, dt) - \lambda \rho(e) de dt,
\end{equation}
where $0<\lambda < \infty$ is the jump intensity and $\rho(e)de$
is the probability measure of each jump size satisfying $\int_E \rho(e)de = 1$. For
jump processes with infinite activities, i.e.,~$\lambda = \infty$,
substantial efforts are needed to construct new spatial discretization approaches,
which is out of scope of this paper and will be considered in our future works.

To proceed, we first introduce a partition the $q$-dimensional Euclidean space $\mathbb{R}^q$ by
$
\mathcal{S} = \mathcal{S}^1 \times \mathcal{S}^2 \times \cdots \times \mathcal{S}^q ,
$
where $\mathcal{S}^k$ for $k = 1, \ldots, q$ is a partition of the one-dimensional space $\mathbb{R}$, i.e.,
\[
\mathcal{S}^k= \Big\{ x_i^k \;\Big|\; x_i^k \in \mathbb{R}, i \in \mathbb{Z}, x_i^k < x^k_{i+1},
\lim_{i\rightarrow +\infty} x_i^k = +\infty, \lim_{i\rightarrow -\infty} x_i^k = -\infty \Big\},
\]
where $\Delta x_k := \max_{i \in \mathbb{Z}}\{|x_i^k - x_{i-1}^k|\}$ and $\Delta x = \max_{1\le k \le q} \Delta x_k$.
For each multi-index $\ii = (i_1, i_2, \ldots, i_q) \in \mathbb{Z}^q$,
the corresponding grid point in $\mathcal{S}$ is denoted by $x_{\ii} = (x^1_{i_1}, \ldots, x^q_{i_q})$.

Recalling \eqref{FK} and \eqref{un}, we can see that, if the terminal condition $\xi$ is a function of $X_T$, $(Y_t^{t,x}, Z_t^{t,x}, \Gamma_{t}^{t, x})$ can be treated as functions of $t$ and $X_t^{t,x} = x$ for $0 \le t \le T$. Analogously, the semi-discrete solution $(Y^n, Z^n, \Gamma^n)$ can be treated as functions of $X^n$. Thus, in this section, we also write $(Y^n, Z^n, \Gamma^n)$ as functions of $x \in \mathbb{R}^q$, i.e.,
\[
\begin{aligned}
Y^n(x) := \mathbb{E}[Y^n | X^n = x],\; Z^n(x) := \mathbb{E}[Z^n | X^n = x],\;  \Gamma^n(x) := \mathbb{E}[\Gamma^n | X^n = x].\\
\end{aligned}
\]
Our objective is to approximate the exact solution $(Y_{t_n}^{t_n,x_\ii}, Z_{t_n}^{t_n,x_\ii}, \Gamma_{t_n}^{t_n, x_\ii})$ by constructing $(Y^n_{\ii}, Z^n_\ii, \Gamma^n_\ii)$, such that
\[
Y^n_{\ii} \approx Y^n(x_\ii) \approx Y_{t_n}^{t_n,x_\ii}, \;\; Z^n_{\ii} \approx Z^n(x_\ii) \approx Z_{t_n}^{t_n,x_\ii},\;\; \Gamma^n_{\ii} \approx \Gamma^n(x_\ii) \approx \Gamma_{t_n}^{t_n,x_\ii},
\]
for $n = 0, \ldots, N-1$ and $\ii \in \mathbb{Z}^q$.

To this end, it is critical to develop effective quadrature rules for approximating the conditional mathematical expectations $\mathbb{E}_{t_n}^{X^n}[\cdot]$ in \eqref{s4:Z}-\eqref{s4:Y}. For instance, at each time-space point $(t_n, x_{\ii}) \in \mathcal{T} \times \mathcal{S}$, approximating $Y^n(x_\ii)$ using \eqref{s4:Y} requires quadrature rules for $\mathbb{E}_{t_n}^{x_\ii}[Y^{n+1}]$ and $\mathbb{E}_{t_n}^{x_\ii}[f^{n+1}]$. In what follows, we take $\mathbb{E}_{t_n}^{x_{\ii}}[Y^{n+1}]$ as an example to propose our new quadrature rule. Slight modifications are needed for approximating $\mathbb{E}_{t_n}^{x_{\ii}}[Y^{n+1}\Delta \tilde{W}^{\top}_{t_{n+1}}]$ and $\mathbb{E}_{t_n}^{x_{\ii}}[Y^{n+1}\Delta \tilde{\mu}^{*}_{t_{n+1}}]$; all the proposed quadrature rules can be directly used to estimate the expectations of $f^{n+1}$, $f^{n+1}\Delta \tilde{W}^{\top}_{t_{n+1}}$ and $f^{n+1}\Delta \tilde{\mu}^{*}_{t_{n+1}}$.

It is observed that $\mathbb{E}_{t_n}^{x_\ii}[Y^{n+1}]$ is defined with respect to the probability measure of the incremental stochastic process
$\Delta X^{n+1} := X^{n+1} - x_{\ii}=  \Phi(t_{n},t_{n+1}, x_{\ii}, I_{\mathcal{J}\in \mathcal A_{\beta}})$ starting from $(t_n,x_{\ii})$, where $\Phi$ is determined by the selected scheme for  the forward SDE. In this section, for the sake of simplicity, we choose the forward Euler method for the scheme in \eqref{s4:X}, i.e.,
\begin{equation}\label{FEU}
\begin{aligned}
X^{n+1}  = x_\ii + b(t_n, x_\ii) \Delta t_n
+ \sigma(t_n, x_\ii) \Delta W_{t_{n+1}}
+ \sum_{k = N_{t_n}+1}^{N_{t_{n+1}}} c(t_n, x_\ii, e_k),\\
\end{aligned}
\end{equation}
where $N_{t}$ for $t \in [t_n, t_{n+1}]$ is the underlying Poisson process.
Since \eqref{FEU} only achieves first-order convergence in the weak sense,
the overall convergence of Scheme 1 will be of first order.
High-order schemes for the forward SDE \cite{Platen:2010eo},
such as order-2.0 weak Taylor scheme, can also be used,
but the corresponding quadrature rules for approximating
$\mathbb{E}_{t_n}^{x_\ii}[\cdot]$ will be dramatically different
from the case of using \eqref{FEU}.
Since the jump intensity $\lambda$ in \eqref{comp} is finite,
the number of jumps of $X_t$ within $(t_n, t_{n+1}]$
follows a compensated Poisson distribution $N_{t_{n+1}} - N_{t_n} - \lambda \Delta t_n$,
where the size of each jump, i.e.,~$c(t_n, x_\ii, e)$,
follows the distribution $\rho(e)de$. Next, we observe that
\begin{equation}\label{dw}
\Delta W^i_{t_{n+1}} = \sqrt{\Delta t_n}\, \xi^i \;\; \text{for}\;\; i = 1, \ldots, d,
\end{equation}
where $\xi^i$ follows the standard normal distribution $N(0,1)$.
Hereafter, we denote by $\varrho(\xi^i)$ the probability density
function of $\xi^i$, and by $\varrho^d(\xi)$ the joint probability density
function of $\xi = (\xi^1, \ldots, \xi^d)^{\top}$.

Now, we can write out the expression of $\mathbb{E}_{t_n}^{x_{\ii}}[Y^{n+1}]$ as
\begin{equation}\label{eq:EY}
\begin{aligned}
&\mathbb{E}_{t_n}^{x_{\ii}} \big[Y^{n+1}\big]\\
= & \sum_{m=0}^\infty \mathbb{P}\Big\{N_{t_{n+1} }- N_{t_n} =m\Big\}\;
\mathbb{E} \bigg[Y^{n+1}\Big(x_{\ii} + b(t_n,x_{\ii})\Delta t_n
+ \sigma(t_n, x_{\ii})\sqrt{\Delta t_n}\, \xi+ \sum_{k=1}^m c(t_n, x_{\ii}, e_k)\Big)\bigg]\\
= & \sum_{m=0}^\infty \mathrm{e}^{-\lambda\Delta t_n} \, \frac{(\lambda \Delta t_n)^m }{m!}\;  \mathbb{E} \bigg[Y^{n+1}\Big(x_{\ii} + b(t_n,x_{\ii})\Delta t_n+\sigma(t_n, x_{\ii}) \sqrt{\Delta t_n}\, \xi + \sum_{k=1}^m c(t_n, x_{\ii}, e_k)\Big)\bigg]\\
= &\, \mathrm{e}^{-\lambda\Delta t_n} \int_{\mathbb{R}^d} Y^{n+1}\Big(x_{\ii} + b(t_n,x_{\ii})\Delta t_n + \sigma(t_n, x_\ii)\sqrt{\Delta t_n}\, \xi\Big)\varrho^d(\xi) d\xi \\
& \hspace{1.5cm} +  \sum_{m=1}^\infty \mathrm{e}^{-\lambda\Delta t_n} \, \frac{(\lambda \Delta t_n)^m }{m!}
\Bigg\{ \int_{\mathbb{R}^d} \int_{E} \cdots \int_{E} Y^{n+1}\Big(x_{\ii} + b(t_n,x_{\ii})\Delta t_n \\
&\hspace{1.5cm}+ \sigma(t_n, x_\ii)\sqrt{\Delta t_n}\, \xi + \sum_{k=1}^m c(t_n, x_\ii,e_k) \Big)
\varrho^d(\xi ) \prod_{k=1}^m \rho(e_k) d\xi\, de_1\cdots de_m\Bigg\},
\end{aligned}
\end{equation}
where $c(t_n, x_\ii, e_k) = c(t_n, x_\ii, e_k^1,\ldots, e_k^q)$
for $k = 1, \ldots, m$ is the size of the $k$-th jump and
$\{c(t_n, x_\ii, e_k)\}_{k=1}^m$ follows the joint distribution
$\prod_{k=1}^m \rho(e_k)$.

Now we study how to approximate $\mathbb{E}_{t_n}^{x_\ii}[Y^{n+1}]$ in \eqref{eq:EY}.
First, we observe that the probability of having $m$ jumps within $(t_n, t_{n+1}]$
is of order $\mathcal{O}((\lambda \Delta t_n)^m)$,
thus the sum of the infinite sequence in \eqref{eq:EY} can be
approximated by the sum of a finite sequence by retaining finite number of jumps.
We denote by ${\mathbb{E}}_{t_n,M_y}^{x_{\ii}}[Y^{n+1}]$ the approximation of $\mathbb{E}_{t_n}^{x_{\ii}}[Y^{n+1}]$ by retaining the first $M_y$ jumps
within $(t_n, t_{n+1}]$. Then, it is easy to see that
the error introduced by the truncation is of order $\mathcal{O}((\lambda \Delta t_n)^{M_y+1})$,
so that $M_y = 2$ is necessary to match the local truncation
error introduced by the semi-discrete scheme in \eqref{s4:Y}.
An analogous notation ${\mathbb{E}}_{t_n,M_f}^{x_{\ii}}[f^{n+1}]$
is used to represent the approximation of ${\mathbb{E}}_{t_n}^{x_{\ii}}[f^{n+1}]$
by retaining the first $M_f$ jumps,
where $M_f = 1$ is sufficient to match the local truncation error in \eqref{s4:Y}.

Next, we also need to approximate a $d$-dimensional integral with respect to
$\xi$ for $m = 0$, and an $m\times q + d$ dimensional
integral with respect to $(\xi,e_1, \ldots, e_m)$ for $m = 1, \ldots, M_y$.
This can be accomplished by selecting an appropriate quadrature rule based on
the properties of $\varrho(\xi)$, $\rho(e)$
and the smoothness of $Y^{n+1}(x)$ with respect to $x$.
A straightforward choice is to use Monte Carlo methods by
drawing samples from $\varrho^d(\xi)$ and $\prod_{k=1}^m \rho(e_k)$,
but they are overall inefficient because of the slow convergence.
When $Y^{n+1}(x)$ is sufficiently smooth with respect to $x$,
an alternative way is to use the tensor product of high-order one-dimensional quadrature rules, e.g.,~Newton-Cotes rules and Gaussian rules, etc.
For example, the integrals with respect to $\xi$ in \eqref{eq:EY}
can be approximated using the tensor product of the Gauss-Hermite rule \cite{Zhao:2010ik}.
For the integrals with respect to $(e_1, \ldots, e_m)$, the Gauss-Legendre rule is a good choice
when $\rho(e)$ is compactly supported, e.g.,~$e$ follows a uniform distribution;
the Gauss-Laguerre rule is appropriate when $\rho(e)$
is the density of an exponential distribution.
Without loss of generality, for $m = 0, \ldots, M_y$,
we denote by $\{w_i^m, s_i^m\}_{i=1}^{S_m}$ and $\{v_j^m, q_j^m\}_{j=1}^{Q_m}$
to represent the chosen quadrature rule for estimating
the integrals in \eqref{eq:EY} with respect to $\xi$ and $(e_1, \ldots, e_m)$,
respectively, where $\{w^m_i\}_{i=1}^{S_m}$, $\{v_j^m\}_{j=1}^{Q_m}$
are quadrature weights and $\{s_i^m\}_{i=1}^{S_m}$, $\{q_j^m\}_{j=1}^{Q_m}$
are quadrature points. Note that $q_j^m$ for $j = 1, \ldots, Q_m$ has $m$ components,
denoted by $\{q_{j,1}^m, \ldots, q_{j,m}^m\}$,
which correspond to the quadrature abscissa for $(e_1, \ldots, e_m)$.
Then, the approximation of $\mathbb{E}_{t_n}^{x_{\ii}}[Y^{n+1}]$,
denoted by $\widehat{\mathbb{E}}_{t_n,M_y}^{x_{\ii}}[Y^{n+1}]$,
is represented by
\vspace{-0.1cm}
\begin{equation}\label{appro_EY}
\begin{aligned}
\widehat{\mathbb{E}}_{t_n,M_y}^{x_{\ii}}\big[Y^{n+1}\big]
\hspace{-0.1cm}:= &\,  \mathrm{e}^{-\lambda \Delta t_n} \sum_{i =1}^{S_0} w_i^0\,
 Y^{n+1}\left( x_{\ii} + b(t_n, x_\ii)\Delta t_n + \sigma(t_n, x_\ii) \sqrt{\Delta t_n} \, s_i^0\right)\\
& + \sum_{m=1}^{M_y} \mathrm{e}^{-\lambda \Delta t_n} \frac{(\lambda \Delta t_n)^m}{m!} \; \Bigg[\sum_{i=1}^{S_m} \sum_{j=1}^{Q_m}
w^m_i v_j^m \;  Y^{n+1}\Big( x_{\ii} + b(t_n, x_\ii)\Delta t_n\\
&\hspace{4cm}  + \sigma(t_n, x_\ii) \sqrt{\Delta t_n}\, s_i^m + \sum_{k=1}^m c(t_n, x_\ii, q_{j,k}^m)\Big)\Bigg].
\end{aligned}
\end{equation}
Analogously, the approximation of $\mathbb{E}_{t_n}^{x_\ii}[f^{n+1}]$,
denoted by $\widehat{\mathbb{E}}_{t_n, M_f}^{x_\ii}[f^{n+1}]$,
can be obtained by replacing $Y^{n+1}$
with $f^{n+1}$ in \eqref{appro_EY}. For the conditional expectations
$\mathbb{E}_{t_n}^{x_\ii}[Y^{n+1}\Delta \tilde{W}_{t_{n+1}}^{\top}]$
and $\mathbb{E}_{t_n}^{x_\ii}[f^{n+1}\Delta \tilde{W}_{t_{n+1}}^{\top}]$,
we observe that each component of
$\Delta \tilde{W}_{t_{n+1}} = (\Delta \tilde{W}_{t_{n+1}}^1,
\ldots, \Delta \tilde{W}_{t_{n+1}}^d)$ defined in \eqref{s3:NBM} can be represented by
\[
\Delta \tilde{W}_{t_{n+1}}^i = \frac{\sqrt{\Delta t_n}}{2} \left(\xi^i + \sqrt{3}\, \tilde{\xi}^i\right) \;\; \text{ for }\; i = 1, \ldots, d,
\]
where $\xi^i$, $\tilde{\xi}^i$ are independent random variables
following standard normal distribution, and $\xi^i$ is the same random variable as
in \eqref{dw}. As such, another $d$-dimensional integral with respect to
$\tilde{\xi} = (\tilde{\xi}^1, \ldots, \tilde{\xi}^d)$ is needed in
\eqref{eq:EY} and \eqref{appro_EY} to define $\widehat{\mathbb{E}}_{t_n, M_y}^{x_\ii}[Y^{n+1}\Delta \tilde{W}_{t_{n+1}}^{\top}]$ and $\widehat{\mathbb{E}}_{t_n, M_f}^{x_\ii}[f^{n+1}\Delta \tilde{W}_{t_{n+1}}^{\top}]$. For the conditional expectations
${\mathbb{E}}_{t_n}^{x_\ii}[Y^{n+1}\Delta \tilde{\mu}^{*}_{t_{n+1}}]$ and ${\mathbb{E}}_{t_n}^{x_\ii}[f^{n+1}\Delta \tilde{\mu}^{*}_{t_{n+1}}]$,
we can see that $\Delta \tilde{\mu}_{t_{n+1}}^{*}$ defined in \eqref{s3:Np} can be represented by
\begin{equation}
\begin{aligned}
\Delta \tilde{\mu}_{t_{n+1}}^* & = \int_{t_n}^{t_{n+1}} \int_E \left(2- \frac{3(t-t_n)}{\Delta t_n}\right)\eta(e)\; \big[{\mu}(de, dt)-\lambda(de)dt\big] \\
& = \sum_{k = N_{t_n}+1}^{N_{t_{n+1}}} \left( 2-\frac{3 (\tau_k - t_n)}{\Delta t_n}\right) \eta(e_k) - \frac{\lambda\Delta t_n}{2} \int_E \eta(e)\rho(e)de,
\end{aligned}
\end{equation}
where $N_t$ for $t \in [t_n, t_{n+1}]$ is the standard Poisson process
and $\tau_k$ for $k = N_{t_n}+1, \ldots, N_{t_{n+1}}$ is the jump time
instant of the $k$-th jump within $(t_n, t_{n+1}]$.
Hence, $\mathbb{E}_{t_n}^{x_\ii}[Y^{n+1} \Delta \tilde{\mu}^{*}_{t_{n+1}}]$
involves another integral with respect to $\tau_k$ compared to
$\mathbb{E}_{t_n}^{x_{\ii}} [Y^{n+1}]$, which requires an additional quadrature rule in \eqref{appro_EY} to construct $\widehat{\mathbb{E}}_{t_n,M_y}^{x_\ii}[{Y}^{n+1}\Delta\tilde{\mu}_{t_{n+1}}^{*}]$.

Based on the quadrature rules used in \eqref{appro_EY},
we observe that it is highly possible the quadrature points
do not belong to the spatial grid $\mathcal{S}$.
In this case, we follow the same strategy as in
\cite{Zhao:2006jx,Zhao:2010ik} to resolve this issue, i.e.,
~constructing piecewise Lagrange interpolating polynomials based on
$\mathcal{S}$ to interpolate the integrands at non-grid quadrature points.
Again, taking $Y^{n+1}(x)$ as an example, it can be approximated by
\[
Y^{n+1}(x) \approx {\widehat{Y}}^{n+1}(x) := \sum_{j_1=1}^{p+1} \cdots \sum_{j_q=1}^{p+1} \Bigg[{Y}_{(i_{j_1}, \ldots, i_{j_q})}^{n+1}
\; \prod_{k = 1}^q \prod_{1\le j \le p+1\atop j\neq j_k} \frac{x^k-x_{i_{j}}^k}{x^k_{i_{j_k}}-x^k_{i_{j}}}\Bigg],
\]
where ${\widehat{Y}}^{n+1}(x)$ is a $p$-th order tensor-product Lagrange interpolating polynomial and $Y^{n+1}_{(i_{j_1}, \ldots, i_{j_q})}$ is the approximate solution of $Y^{n+1}(x)$ at the spatial point $(x_{i_{j_1}}^1, \ldots, x^q_{i_{j_q}})$.
For $k = 1,\ldots, q$, the interpolation points
$\{x^k_{i_{j}}\}_{j=1}^{p+1} \subset \mathcal{S}^k$ are the {\em closest}
$p+1$ neighboring points of $x^k$, such that $(x^1_{i_{j_1}}, \ldots, x^q_{i_{j_q}})$
for $j_k = 1, \ldots, p+1$ and $k = 1, \ldots, q$ constitute a local tensor-product
sub-grid around $x$.
In summary, the fully discrete scheme of the FBSDEs in \eqref{s1:e1} is given as follows:
\begin{scheme}\label{s4:sch:2}
Given initial condition $X_0$ for the forward SDE in \eqref{s1:e1}
and the terminal condition $\varphi(X_T)$ for the backward SDE in \eqref{s1:e1},
solve the approximate solution $(Y_{\ii}^n, Z_\ii^n, \Gamma_\ii^n)$, for $n=N-1,\cdots,0$ and $\ii \in \mathbb{Z}^q$, by
\begin{align}
\label{s4:X1}
X^{n+1} & =  x_\ii + \Phi(t_{n},t_{n+1}, x_\ii, I_{\mathcal{J}\in \mathcal A_{\beta}}),\\
\label{s4:Y1}
Y^n_\ii & = \widehat{\mathbb{E}}_{t_n,M_y}^{x_\ii}\left[\widehat{Y}^{n+1}\right]
+\frac{1}{2}\Delta t_n f^n_\ii
+ \frac{1}{2}\Delta t_n \widehat{\mathbb{E}}_{t_n,M_f}^{x_\ii}\left[\widehat{f}^{n+1}\right],\\
\label{s4:Z1}
\frac{1}{2}\Delta t_n Z^n_\ii
& =  \widehat{\mathbb{E}}_{t_n,M_y}^{x_\ii}\left[\widehat{Y}^{n+1}\Delta \tilde{W}_{t_{n+1}}^{\top}\right]
+\Delta t_n \widehat{\mathbb{E}}_{t_n,M_f}^{x_\ii}\left[\widehat{f}^{n+1}\Delta \tilde{W}_{t_{n+1}}^{\top}\right],\\
\label{s4:G1}
\frac{1}{2}\Delta t_n \Gamma^n_\ii
& = \widehat{\mathbb{E}}_{t_n,M_y}^{x_\ii}\left[\widehat{Y}^{n+1}
\Delta\tilde{\mu}_{t_{n+1}}^{*}\right]
+\Delta t_n
\widehat{\mathbb{E}}_{t_n,M_f}^{x_\ii}
\left[\widehat{f}^{n+1}\Delta\tilde{\mu}_{t_{n+1}}^{*}\right],
\end{align}
where $\widehat{f}^{n+1} = f(t_{n+1},X^{n+1}, \widehat{Y}^{n+1}, \widehat{Z}^{n+1},
\widehat{\Gamma}^{n+1})$, $f_{\ii}^n = (t_n, x_\ii, Y^n_\ii, Z^{n}_\ii, \Gamma^{n}_\ii)$,
$\Delta \tilde W_{t_{n+1}}^{\top}$ and $\Delta\tilde{\mu}_{t_{n+1}}^{*}$
are defined in \eqref{s3:NBM} and \eqref{s3:Np} with $s = t_{n+1}$, respectively.
\end{scheme}

Similar to the semi-discrete scheme, Scheme \ref{s4:sch:2} can be directly used as a fully discrete scheme for the PIDE in \eqref{eq:Kol}.
The solution $u(t_n,x_{\ii})$ of the PIDE is approximated by $Y_{\ii}^n$ for $n = 0, \ldots, N-1$ and $\ii \in \mathbb{Z}^q$.
We observe that at each grid point $(t_n, x_\ii)$, the computation of $Y_\ii^n$
only depends on $(X^{n+1},Y^{n+1},Z^{n+1},\Gamma^{n+1})$ even though an
implicit time-stepping scheme is used. This means $\{Y_{\ii}^n\}_{\ii \in \mathbb{Z}^q}$
at each time step can be computed {\em independently},
so that the difficulty of solving linear systems with possibly {\em dense} matrices,
due to the nonlocality of the integral operator, is completely avoided.
This feature makes it straightforward to develop massively parallel
algorithms and incorporate adaptive spatial interpolation methods.
\begin{rem}
It is noted that the total computational cost of the Scheme \ref{s4:sch:2}
is dominated by the cost of approximating
$\mathbb{E}_{t_n}^{x_\ii}[\cdot]$ at each grid point
$(t_n,x_{\ii}) \in \mathcal{T}\times \mathcal{S}$ using the formula in \eqref{appro_EY}.
For example, when solving a three-dimensional problem $q=d=3$
and retaining two L\`{e}vy jumps $M_y = M_f =2$, we are
facing a large amount of six-dimensional integration problems.
In this case, sparse-grid quadrature rules
\cite{Bungartz:2004kx,Gunzburger:2014hi,Zhang:2013en}
can be used to alleviate the explosion of computational cost due to {\em curse of dimensionality}.
\end{rem}
\vspace{-0.2cm}
\section{Numerical examples}\label{sec:num}
\vspace{-0.2cm}
In this section, we report on the results of two one-dimensional
numerical examples that illustrate the accuracy
and the effectiveness of Schemes \ref{s4:sch:1} and \ref{s4:sch:2}.
We take uniform partitions in both temporal and
spatial domains with the time and space step sizes denoted by
$\Delta t$ and $\Delta x$, respectively.
The time step number $N$ is then given by $N = {T}/{\Delta t}$ where $T$ is the terminal time.
For the sake of illustration, we only solve FBSDEs on bounded spatial domains.
%Piecewise Lagrange interpolation and
%tensor-product quadrature rules are used to approximate
%$\mathbb{E}_{t_n}^{x_i}[\cdot]$, where the one-dimensional
%quadrature rule is chosen according to the property of
%the L\`{e}vy measure $\lambda(de)$, the terminal condition $\varphi$
%and the generator $f$.
The goal is to test the convergence
rates of time discertizaiton and spatial interpolation with respect to
$\Delta t$ and $\Delta x$, respectively.
To this end, we always set the number of quadrature points
to be sufficiently large, so that the error contributed by
the use of quadrature rules is too small to affect the convergence rates of interest.

\vspace{-0.4cm}
\subsection{Example 1}\label{sec:ex1}
\vspace{-0.3cm}
We consider the following {\em nonlinear} FBSDEs:
\begin{equation}\label{ex1:eq1}
\left\{
\begin{aligned}
 d X_t & = dW_t + \int_E e \, \tilde{\mu}(de, dt), \\
 -dY_t & = \Bigg\{\dfrac{(Y_t-2) \cdot\exp(Y_t)}{2\, \exp\big[\sin(X_t+t)+2\big]}
 - \frac{Z_t\, Y_t}{\sin(X_t+t)+2}  - \Gamma_t\Bigg\} dt
          - Z_t \, dW_t - \int_E U(e) \tilde{\mu}(de,dt),
\end{aligned}\right.
\end{equation}
where the terminal condition is $\varphi(X_T) = \sin(X_T + T)+2$.
The L\`{e}vy measure is defined by
\begin{equation}\label{levym}
\lambda(de) = \lambda \rho(e) de
:= \mathcal{X}_{[-\delta,\delta]} (e) de \; \mbox{ with }\; \delta >0,
\end{equation}
where $\mathcal{X}_{[-\delta, \delta]}(e)$
is the characteristic function of the interval $[-\delta, \delta]$,
so that $\lambda = 2 \delta$ is the jump intensity and
$\rho(e) = \frac{1}{2\delta} \mathcal{X}_{[-\delta, \delta]}(e)$
is the density function of a uniform distribution on $[-\delta, \delta]$.
The exact solution of the FBSDEs is
\begin{equation}
\left\{\begin{aligned}
Y_t & = \sin(X_t + t) + 2,\\
Z_t & = \cos(X_t + t),\\
\Gamma_t & = \cos(X_t + t -\delta) - \cos(X_t + t + \delta) - 2 \delta \sin(X_t + t).
\end{aligned}\right.
\end{equation}
Accordingly, the PIDE corresponding to \eqref{ex1:eq1} is
\[
\left\{
\begin{aligned}
& \frac{\partial u}{\partial t} + \frac{1}{2} \frac{\partial^2 u}{\partial x^2}+ \int_{E} \left(u(t,x + e) -u(t,x)
 \right) \lambda (de) \\
 &\hspace{0.5cm}+ \dfrac{(u-2) \cdot\exp(u)}{2\, \exp\big[\sin(x+t)+2\big]}- \frac{u}{\sin(x+t)+2}\dfrac{\partial u}{\partial x}  - \mathcal{B}[u] = 0,\\
& u(T,x) = \sin(x+T) + 2,
\end{aligned}\right.
\]
where $\mathcal{B}[u] = \cos(x + t -\delta) - \cos(x + t + \delta) - 2 \delta \sin(x + t)$.

Since the density function $\rho(e)$ is uniform with
the support $[-\delta, \delta]$, we use the tensor product
of the 8-point Gauss-Legendre rule and the 8-point
Gauss-Hermite rule to approximate the integrals involved
in $\mathbb{E}_{t_n}^{x_i}[\cdot]$.
%We observe that the forward SDE
%is just the Brownian motion plus the standard compensated
%compound Poisson process where the diffusion coefficient $\sigma$
%and the jump size $c$ are independent of $t$ and $X_t$.
%As such, the discretization of the forward SDE is not needed.

First, we test the convergence rate with respect to $\Delta t$
where the terminal time is $T=1$. To this end, we set $\Delta x = 0.01$
and use piecewise cubic Lagrange interpolation to construct $\widehat{Y}^{n+1}(x)$
for $n = 0, \ldots, N-1$, such that the time discretization
error dominates the total error. Setting $\delta = 1$ and
$\Delta t = 2^{-4}, 2^{-5}, 2^{-6}, 2^{-7}, 2^{-8}$,
the numerical results are shown in Table
\ref{ex1:t1}. As expected, the convergence rate with respect to
$\Delta t$ depends on the number of jumps retained in
$\widehat{\mathbb{E}}_{t_n,M_y}^{x_\ii}[\cdot]$ and
$\widehat{\mathbb{E}}_{t_n,M_f}^{x_\ii}[\cdot]$.
For example, when $M_y = M_f = 0$, i.e., no jump is included,
our scheme fails to converge. In order to achieve second-order
convergence, we must set $M_y \ge 2$ and $M_f \ge 1$.
\begin{table}[h!]
 \scriptsize
\begin{center}
\caption{Errors and convergence rates with respect to $\Delta t$ in Example 1,
where $T =1$, $\delta = 1$, $\Delta x = 0.01$,
and piecewise cubic Lagrange interpolation is used.} \label{ex1:t1}
\begin{tabular}{|c|c|c|c|c|c|c|}\hline
\multicolumn{7}{|c|}{$\| Y_{0}^{0,x} -\widehat{Y}^{0}(x)\|_{L^{\infty}([0,1])}$}\\
\hline   & $\Delta t = 2^{-4}$ & $\Delta t = 2^{-5}$ & $\Delta t = 2^{-6}$ & $\Delta t = 2^{-7}$ & $\Delta t = 2^{-8}$ & CR \\
 \hline
 $M_y = 0$, $M_f=0$ &  5.178E-1  &   4.786E-1  &   4.498E-1  &   4.290E-1  &   4.090E-1  &   0.084  \\
  \hline
 $M_y = 1$, $M_f=0$  &  4.155E-2  &   1.778E-2  &   7.855E-3  &   3.631E-3  &   1.691E-3  &   1.153  \\
 \hline
 $M_y = 2$, $M_f=1$ &  4.539E-3  &   9.878E-4  &   2.211E-4  &   5.065E-5  &   1.144E-5  &   2.155  \\
\hline
 $M_y = 3$, $M_f=2$ &  3.414E-3  &   7.305E-4  &   1.609E-4  &   3.646E-5  &   8.087E-6  &   2.177  \\
\hline
\hline
\multicolumn{7}{|c|}{$\| Z_{0}^{0,x} -\widehat{Z}^{0}(x)\|_{L^{\infty}([0,1])}$}\\
\hline   & $\Delta t = 2^{-4}$ & $\Delta t = 2^{-5}$ & $\Delta t = 2^{-6}$ & $\Delta t = 2^{-7}$ & $\Delta t = 2^{-8}$ & CR \\
 \hline
 $M_y = 0$, $M_f=0$ &  1.591E-0  &   2.155E-0  &   2.505E-0  &   2.237E-0  &   2.534E-0  &   -0.140  \\
  \hline
 $M_y = 1$, $M_f=0$  &  1.629E-1  &   9.244E-2  &   5.028E-2  &  2.156E-2  &   1.146E-2  &   0.976  \\
 \hline
 $M_y = 2$, $M_f=1$ &  1.846E-2  &   5.459E-3  &   1.536E-3  &   3.293E-4  &   8.706E-5  &   1.951  \\
\hline
 $M_y = 3$, $M_f=2$ &  1.475E-2  &   4.230E-3  &  1.161E-3  &   2.450E-4  &   6.376E-5  &   1.982  \\
\hline
\hline
\multicolumn{7}{|c|}{$\| \Gamma_{0}^{0,x} -\widehat{\Gamma}^{0}(x)\|_{L^{\infty}([0,1])}$}\\
\hline   & $\Delta t = 2^{-4}$ & $\Delta t = 2^{-5}$ & $\Delta t = 2^{-6}$ & $\Delta t = 2^{-7}$ & $\Delta t = 2^{-8}$ & CR \\
 \hline
 $M_y = 0$, $M_f=0$ &  5.165E-1  &   7.206E-1  &   6.360E-1  &   5.444E-1  &   5.270E-1  &   0.035  \\
  \hline
 $M_y = 1$, $M_f=0$  &   3.519E-1  &   1.805E-1  &   9.636E-2  &   4.416E-2  &   2.240E-2  &   0.998 \\
 \hline
 $M_y = 2$, $M_f=1$ &  2.151E-2  &   5.364E-3  &   1.453E-3  &   3.365E-4  &   8.458E-5  &   1.998  \\
\hline
 $M_y = 3$, $M_f=2$ &  1.549E-2  &   3.856E-3  &   1.057E-3  &   2.467E-4  &   6.212E-5  &   1.989  \\
\hline
\end{tabular}
\end{center}
\end{table}
\begin{table}[h!]
 \scriptsize
\begin{center}
\caption{Errors and convergence rates with respect to $\Delta x$ in Example 1,
where $T =1$, $\delta = 1$, $x \in [0,1]$, $N = 1024$, $M_y = 3$ and $M_f = 2$.} \label{ex1:t3}
\begin{tabular}{|c|c|c|c|c|c|c|}\hline
\multicolumn{7}{|c|}{Linear interpolation}\\
\hline   & $\Delta x = 2^{-2}$ & $\Delta x = 2^{-3}$ & $\Delta x = 2^{-4}$ & $\Delta x = 2^{-5}$ & $\Delta x = 2^{-6}$  & CR \\
 \hline
$\| Y_{0}^{0,x} -\widehat{Y}^{0}(x)\|_{\infty}$ & 1.756E-2 & 6.094E-3  &   1.968E-3  &  3.194E-4  &   1.012E-4  &     2.036  \\
  \hline
$\| Z_{0}^{0,x} -\widehat{Z}^{0}(x)\|_{\infty}$ & 9.273E-2 &  3.163E-2  &   8.457E-3  &   2.278E-3  &   5.203E-4  &    1.885   \\
\hline
$\| \Gamma_{0}^{0,x} -\widehat{\Gamma}^{0}(x)\|_{\infty}$ & 2.000E-2 & 6.812E-3  &   2.173E-3  &   3.360E-4  &   1.081E-4  &   2.063 \\
\hline
\hline
\multicolumn{7}{|c|}{Quadratic interpolation}\\
\hline   & $\Delta x = 2^{-2}$ & $\Delta x = 2^{-3}$ & $\Delta x = 2^{-4}$ & $\Delta x = 2^{-5}$ & $\Delta x = 2^{-6}$  & CR \\
 \hline
$\| Y_{0}^{0,x} -\widehat{Y}^{0}(x)\|_{\infty}$  &  6.183E-2  &   8.004E-3  &   9.941E-4  &   1.266E-4  & 1.438E-5  &     3.012  \\
  \hline
$\| Z_{0}^{0,x} -\widehat{Z}^{0}(x)\|_{\infty}$ &  5.537E-2  &   7.016E-3  &   1.583E-3  &   1.266E-4  & 1.096E-5  &  3.039   \\
\hline
$\| \Gamma_{0}^{0,x} -\widehat{\Gamma}^{0}(x)\|_{\infty}$ & 1.927E-2  &   3.489E-3  &   3.527E-4  &   5.531E-5 & 5.587E-6 &   2.948 \\
\hline
\end{tabular}
\end{center}
\end{table}

Next, we test the convergence rate with respect to $\Delta x$ by setting $\delta = 1$, $T = 1$, $N = 1024$, $M_y = 3$, $M_f = 2$, and $\Delta x = 2^{-2}, 2^{-3}, 2^{-4}, 2^{-5}, 2^{-6}$. The error is measured in $L^{\infty}$ norm. In Table \ref{ex1:t3}, we can see that the spatial discretization error decays as expected, i.e.,~second-order and third-order convergence rates for piecewise linear and piecewise quadratic interpolations, respectively.

\subsection{Example 2}\label{sec:ex2}
 We consider the following nonlinear FBSDE:
\begin{equation}\label{ex2:eq}
\left\{
\begin{aligned}
 d X_t & = \sin(2X_t+t) dt + \big[\cos(X_t)+t+2\big]dW_t + \int_E e \, \tilde{\mu}(de, dt), \\
 -dY_t & = \Bigg\{ -\frac{\sin(2X_t+t)Y_tZ_t}{[\cos(X_t)+t+2](\sin(t)+2)\exp(-X_t)} \\
          &\hspace{0.7cm} - 0.5[\cos(X_t)+t+2]^2 Y_t- \Gamma_t\Bigg\} dt  - Z_t \, dW_t - \int_E U(e) \tilde{\mu}(de,dt),
\end{aligned}\right.
\end{equation}
where the terminal condition is $\varphi(X_T) = [\sin(T)+2] \exp(-X_T)$. The L\`{e}vy measure $\tilde{\mu}(de,dt)$ is defined as in
\eqref{levym}. The exact solution of the FBSDEs is
\begin{equation}
\left\{\begin{aligned}
Y_t & = (\sin(t)+ 2) \exp(-X_t),\\
Z_t & = -(\cos(X_t)+t+2)(\sin(t) + 2) \exp(-X_t),\\
\Gamma_t & = (\sin(t)+2)\big[\exp(-X_t +\delta) - \exp(-X_t-\delta) - 2\delta\exp(-X_t)\big].
\end{aligned}\right.
\end{equation}
Accordingly, the PIDE corresponding to \eqref{ex2:eq} is
\[
\left\{
\begin{aligned}
& \frac{\partial u}{\partial t} +  \sin(2x+t)\frac{\partial u}{\partial x} + \frac{1}{2} \big[\cos(x)+t+2\big]^2\frac{\partial^2 u}{\partial x^2}\\
&\hspace{0.5cm} + \int_{E} \left(u(t,x + e) -u(t,x)
 \right) \lambda (de)  -\frac{\sin(2x+t)u(t,x) }{(\sin(t)+2)\exp(-x)} \frac{\partial u}{\partial x} \\
 &\hspace{0.5cm} -  \frac{1}{2}[\cos(x)+t+2]^2 u(t,x)- \mathcal{B}[u] = 0,\\
& u(T,x) = \sin(T+2) \exp(-x),
\end{aligned}\right.
\]
where $\mathcal{B}[u] = (\sin(t) + 2)[\exp(-x+\delta)-\exp(-x-\delta)-2\delta \exp(-x)]$.
Similar to Example 1,
we use the tensor product of the 8-point Gauss-Legendre rule
and the 8-point Gauss-Hermite rule to approximate the integrals involved in $\mathbb{E}_{t_n}^{x_i}[\cdot]$.
%Observing that both the drift and diffusion coefficients
%of the forward SDE in \eqref{ex2:eq} depend on $t$ and $X_t$,
%a numerical scheme is needed for discretizing the forward SDE in \eqref{ex2:eq}.
In this example, we use forward Euler scheme in \eqref{FEU},
such that the overall convergence rate will be expected to be first order.

First, we test the convergence rate with respect to $\Delta t$ where the terminal time is $T=1$. To this end, we set $\Delta x = 0.01$ and use piecewise cubic Lagrange interpolation to construct $\widehat{Y}^{n+1}(x)$ for $n = 0, \ldots, N-1$, so that the time discretization error dominates the total error. Setting $\delta = 1$ and $\Delta t = 2^{-5}, 2^{-6}, 2^{-7}, 2^{-8}, 2^{-9}$, the numerical results are shown in Table
\ref{ex2:t1}. As expected, the convergence rate with respect to $\Delta t$ depends on the number of jumps retained in constructing $\widehat{\mathbb{E}}_{t_n,M_y}^{x_\ii}[\cdot]$ and $\widehat{\mathbb{E}}_{t_n,M_f}^{x_\ii}[\cdot]$. In this case, we can only achieve, at most,
first-order convergence with respect to $\Delta t$ due to the use of the forward Euler scheme.

Next, we test the convergence rate with respect to $\Delta x$ by setting $\delta = 1$, $T = 1$, $N = 1024$, $M_y = 2$ and $M_f =1$. The spatial mesh size is set to $\Delta x = 2^{-5}, 2^{-6}, 2^{-7}, 2^{-8}, 2^{-9}$ for linear interpolation and
$\Delta x = 2^{-2}, 2^{-3}, 2^{-4}, 2^{-5}, 2^{-6}$ for quadratic interpolation. The error is measured in $L^{\infty}$ norm. In Table \ref{ex2:t3}, we can see that the spatial discretization error decays as expected, i.e.,~second-order and third-order convergence rates for piecewise linear and piecewise quadratic interpolations, respectively.

\begin{table}[h!]
 \scriptsize
\begin{center}
\caption{Errors and convergence rates with respect to $\Delta t$ in Example 2,
where $T =1$, $\delta = 1$, $\Delta x = 0.01$,
and piecewise cubic Lagrange interpolation are used.} \label{ex2:t1}
\begin{tabular}{|c|c|c|c|c|c|c|}\hline
\multicolumn{7}{|c|}{$\| Y_{0}^{0,x} -\widehat{Y}^{0}(x)\|_{L^{\infty}([0,1])}$}\\
\hline   &  $\Delta t = 2^{-5}$ & $\Delta t = 2^{-6}$ & $\Delta t = 2^{-7}$ & $\Delta t = 2^{-8}$ & $\Delta t = 2^{-9}$ &CR \\
 \hline
 $M_y = 0$, $M_f=0$ &  1.331E-1  &   9.191E-2  &   7.689E-2  &   6.763E-2  &   6.206E-2  &   0.264  \\
  \hline
 $M_y = 1$, $M_f=0$  &  2.996E-2  &   1.362E-2  &   5.339E-3  &   2.226E-3  &   1.998E-3  &   1.043 \\
 \hline
 $M_y = 2$, $M_f=1$ &  3.071E-2  &   1.047E-2  &   3.835E-3  &   1.613E-3  &  7.155E-4  &  1.355  \\
\hline
\hline
\multicolumn{7}{|c|}{$\| Z_{0}^{0,x} -\widehat{Z}^{0}(x)\|_{L^{\infty}([0,1])}$}\\
\hline   &  $\Delta t = 2^{-5}$ & $\Delta t = 2^{-6}$ & $\Delta t = 2^{-7}$ & $\Delta t = 2^{-8}$ & $\Delta t = 2^{-9}$ &CR \\
 \hline
 $M_y = 0$, $M_f=0$ &  4.757E-1  &   4.522E-1  &   5.558E-1  &   6.583E-1  &   6.881E-1 &   -0.161  \\
  \hline
 $M_y = 1$, $M_f=0$  &  3.779E-1  &   1.733E-1  &   7.756E-2  &  3.598E-2  &   1.661E-2  &   1.128  \\
 \hline
 $M_y = 2$, $M_f=1$ &  1.129E-1  &   5.351E-2  &   2.714E-2  &   1.201E-2  &   5.771E-3  &   1.074  \\
\hline
\hline
\multicolumn{7}{|c|}{$\| \Gamma_{0}^{0,x} -\widehat{\Gamma}^{0}(x)\|_{L^{\infty}([0,1])}$}\\
\hline   &  $\Delta t = 2^{-5}$ & $\Delta t = 2^{-6}$ & $\Delta t = 2^{-7}$ & $\Delta t = 2^{-8}$ & $\Delta t = 2^{-9}$ &CR \\
 \hline
 $M_y = 0$, $M_f=0$ &  1.379E-1  &   1.194E-1  &   8.948E-2  &   8.391E-2  &   7.765E-2  &  0.217  \\
  \hline
 $M_y = 1$, $M_f=0$  &   5.125E-2  &   2.528E-2  &   1.226E-2  &   7.598E-3  &   3.996E-3  &   0.909 \\
 \hline
 $M_y = 2$, $M_f=1$ &  4.789E-2  &   2.116E-2  &   1.057E-2  &   5.373E-3  &   2.509E-3  &   1.049  \\
\hline
\end{tabular}
\end{center}
\end{table}

\begin{table}[h!]
 \scriptsize
\begin{center}
\caption{Errors and convergence rates with respect to $\Delta x$ in Example 2,
where $T =1$, $\delta = 1$, $x\in [0,1]$, $N = 1024$, $M_y = 2$ and $M_f = 1$.} \label{ex2:t3}
\begin{tabular}{|c|c|c|c|c|c|c|}\hline
\multicolumn{7}{|c|}{Linear interpolation}\\
\hline   & $\Delta x = 2^{-5}$ & $\Delta x = 2^{-6}$ & $\Delta x = 2^{-7}$ & $\Delta x = 2^{-8}$ & $\Delta x = 2^{-9}$  & CR \\
 \hline
$\| Y_{0}^{0,x} -\widehat{Y}^{0}(x)\|_{\infty}$ & 2.981E-3 & 8.383E-4  &   2.227E-4  &  5.583E-5  &   1.652E-5  &     1.899  \\
  \hline
$\| Z_{0}^{0,x} -\widehat{Z}^{0}(x)\|_{\infty}$ & 8.538E-2 &  2.774E-2  &   6.109E-3  &   1.810E-3  &   4.206E-4  &    1.926   \\
\hline
$\| \Gamma_{0}^{0,x} -\widehat{\Gamma}^{0}(x)\|_{\infty}$ & 5.859E-3 & 1.537E-3  &   2.667E-4  &   7.729E-5  &   2.097E-5  &   2.056 \\
\hline
\hline
\multicolumn{7}{|c|}{Quadratic interpolation}\\
\hline   & $\Delta x = 2^{-2}$ & $\Delta x = 2^{-3}$ & $\Delta x = 2^{-4}$ & $\Delta x = 2^{-5}$ & $\Delta x = 2^{-6}$  & CR \\
 \hline
$\| Y_{0}^{0,x} -\widehat{Y}^{0}(x)\|_{\infty}$ &  3.127E-2  &   3.916E-3  &   6.021E-4  &   8.410E-5  & 1.093E-5  &   2.856  \\
  \hline
$\| Z_{0}^{0,x} -\widehat{Z}^{0}(x)\|_{\infty}$ &  8.945E-2  &   1.158E-2  &   1.376E-3  &   1.966E-4  & 2.638E-5  &  2.935   \\
\hline
$\| \Gamma_{0}^{0,x} -\widehat{\Gamma}^{0}(x)\|_{\infty}$ & 7.512E-3  &  1.349E-3  &  2.109E-4  &   3.236E-5 & 2.807E-6 &   2.815\\
\hline
\end{tabular}
\end{center}
\end{table}

\section{Concluding remarks}\label{sec:con}
In this work, we propose new numerical schemes for decoupled forward-backward stochastic differential equations with jumps, which feature high-order temporal and spatial convergence rates. This advantage has been verified by both theoretical analysis and numerical experiments. Meanwhile, we also realized that our schemes cannot achieve the desired convergence rates in the sense that
the solution of the FBSDEs does not satisfy the necessary regularity conditions. For example, this may happen in real-world financial problems, such as option pricing. However, the regularity conditions do not limit the applicability of the proposed approach, because our method can be directly employed as a probabilistic scheme for related PIDEs which are widely used to describe anomalous transport in subsurface flow and plasma physics. In these settings, there is a variety of problems satisfying the regularity conditions, and high-order schemes are highly desired. Moreover, compared to existing deterministic approaches (e.g.,~finite elements) for the PIDEs, the ability to completely avoids the solution of {\em dense} linear systems, as well as to utilize efficient adaptive approximation, and the potential of massively parallel implementation, make our technique highly advantageous. Our future works will focus on extending the proposed numerical schemes to the case of Poisson random measures with {\em infinite} activities, and integrating sparse grid methods for high-dimensional FBSDEs with jumps.\\

\noindent{\bf Acknowledgments.} This work is partially supported by the
National Natural Science Foundations of China under grant numbers 91130003
and 11171189; by Natural Science Foundation of
Shandong Province under grant number ZR2011AZ002;
by the U.S.~Department of Energy, Office of Science,
Office of Advanced Scientific Computing Research,
Applied Mathematics program under contract number ERKJE45;
and by the Laboratory Directed Research and Development program
at the Oak Ridge National Laboratory, which is operated
by UT-Battelle, LLC, for the U.S.~Department of Energy
under Contract DE-AC05-00OR22725.

\end{document}